\documentclass[12pt,a4paper]{article}
\usepackage[]{inputenc}
\usepackage{amsmath}
\usepackage{float}
\usepackage{multirow}
\usepackage{array}
\usepackage{amsfonts}
\usepackage{amssymb}
\usepackage{color}
\usepackage{authblk}
\usepackage{graphicx}
\usepackage[labelfont=bf]{caption}
\usepackage{caption,subcaption}
\usepackage{hyperref}
\usepackage[left=2cm,right=2cm,top=2cm,bottom=2cm]{geometry}

\graphicspath{{figures/}}

\begin{document}
\title{Stabilized reduced basis methods for parametrized steady Stokes and Navier-Stokes equations}

\author[1, 2]{\small Shafqat Ali\thanks{ali.qau1987@gmail.com, shafqat.ali@giki.edu.pk}}
\author[1]{\small Francesco Ballarin\thanks{francesco.ballarin@sissa.it}}
\author[1]{\small Gianluigi Rozza\thanks{gianluigi.rozza@sissa.it (corresponding author)}}
\affil[1]{\footnotesize mathLab, Mathematics area, SISSA, via Bonomea 265, I-34136, Trieste, Italy}
\affil[2]{\footnotesize Current address: Faculty of Engineering Sciences, Ghulam Ishaq Khan Institute of Engineering Sciences and Technology, Topi, 44000 Pakistan}
\date{}
\maketitle
\begin{abstract}
It is well known in the Reduced Basis approximation of saddle point problems that the Galerkin projection on the reduced space does not guarantee the inf-sup approximation stability even if a stable high fidelity method was used to generate snapshots. For problems in computational fluid dynamics, the lack of inf-sup stability is reflected by the inability to accurately approximate the pressure field. In this context, inf-sup stability is usually recovered through the enrichment of the velocity space with suitable supremizer functions. The main goal of this work is to propose an alternative approach, which relies on the residual based stabilization techniques customarily employed in the Finite Element literature, such as Brezzi-Pitkaranta, Franca-Hughes, streamline upwind Petrov-Galerkin, Galerkin Least Square.
In the spirit of \textit{offline-online} reduced basis computational splitting, two such options are proposed, namely \textit{offline-only stabilization} and \textit{offline-online stabilization}. These approaches are then compared to (and combined with) the state of the art supremizer enrichment approach. Numerical results are discussed, highlighting that the proposed methodology allows to obtain smaller reduced basis spaces (i.e., neglecting supremizer enrichment) for which a modified inf-sup stability is still preserved at the reduced order level.
\end{abstract}  


\section{Introduction}
\label{sec:introduction}
Numerical approximation of fluid dynamics problems (modelled by incompressible Stokes and Navier-Stokes equations) is required by several applications in mechanics and engineering, often depending on some parameters. In case of repeated evaluation for different parametric configurations, the Reduced Basis (RB) method \cite{RB2016} offers attractive performance, able to cut down computational costs of standard Finite Element (FE) \cite{QV} simulations. In the context of fluid dynamics problems, formulations based on a Galerkin projection, such as FE and RB, require the fulfillment of a discrete $\inf$-$\sup$ condition \cite{boffi2013mixed}. While several choices of suitable discrete spaces for the velocity and pressure are well known at the FE level \cite{boffi2013mixed,QV}, the corresponding RB formulation still deserves investigation, even after more than a decade from the original supremizer enrichment proposal \cite{Rozza2005a,Rovas2003,Veroy2007}.
This approach consists in the introduction of the inner pressure supremizer for the velocity-pressure stability of the RB spaces.
Several works on RB method for Stokes problem \cite{NegriManzoniRozza2015, Rozza2013} were developed subsequently using the pressure stabilization via the inner pressure supremizer operator. This approach was then extended to Navier-Stokes problems as well \cite{Simone2008,Simone2009,Manzoni2014,Veroy2005, Ballarin2015}. \par 
However, even though nowadays the supremizer enrichment approach is widely used in the RB community, it still carries the side effect of increasing the RB velocity space, resulting in a possible increase in the online computational time. Therefore, in this work we look for an alternative approach which does not entail an enrichment procedure, yet results in a stable RB approximation. 
\par
The stabilization approach developed in this work builds upon classical residual based stabilization techniques \cite{Hughes1986, Brooks1982, Hughes1989, Hughes1980, Douglas1989, ZHU1993} already known for FE methods, and properly adapted to a reduced order setting in this work. Such stabilization techniques have been employed to handle two sources of instabilities. The first one (which motivates our work) is the incompatibility of velocity and pressure FE pairs. Another possible instability is due to dominating advection for problems characterized by large Reynolds numbers. In this work we will only focus on the inf-sup stability issue, referring to \cite{PR014a, Traian2015, Codina2013, Balajewicz2012, Torlo2018} for reduced order methods for advection dominated problems.
\par
The idea of combining stabilization techniques to model reduction methods follows the work of Pacciarini and Rozza \cite{PR014a} and Torlo et al. \cite{Torlo2018} for advection-dominated problem, where they introduced the concept of \textit{offline-online stabilization} and the \textit{offline-only stabilization}. 
The main novelty here is the comparison and combination to the supremizer approach by Rozza et al. \cite{Veroy2007}, as both methodologies aim at preserving the $\inf$-$\sup$ stability of the reduced order system.
As in \cite{PR014a,PP015,PR014b}, the \textit{offline-online stabilization} method is based on performing the Galerkin projection with respect to the stabilized formulation in both \textit{offline} and \textit{online} stages. In contrast, the \textit{offline-only stabilization} consists in using the stabilized formulation only during the \textit{offline} stage, while projecting with respect to the standard formulation during the \textit{online} stage. 
\par
We refer to some related works in recent past on the stabilization of reduced order models, mostly based on POD, see for instance Baiges et al. \cite{Baiges2014, Codina2013}, Caiazzo et al. \cite{Traian2014} and L{\o}vgren et al. \cite{Deparis2012}. The difference in our approach with respect to previous works is that the stabilization approach we are using is residual based and is a strongly consistent stabilization. 
\par 
This manuscript is mainly divided into two parts, with the focus on steady Stokes problem in first part and steady Navier-Stokes problem in second part, respectively. After this introduction, in section \ref{sec: S_continuous_formulation}, we present FE and RB approximation for the steady Stokes problem in a parametrized domain, with a particular focus on the introduction of residual based stabilization techniques in the formulation. In section \ref{sec: Numerical_Results} we show some numerical results and error analysis to compare different stabilization options.
In the second part, section \ref{sec: NS_continuous_formulation} shows FE and RB discretizations for the steady Navier-Stokes problem. Afterwards, in section \ref{sec:Results_NS} we show some numerical results on a benchmark case, addressing the performance (in terms of cost and accuracy) of different stabilization approaches presented. Finally some conclusions and perspectives are provided in section \ref{sec: Conclusion}.
\section{Parametrized steady Stokes problem}
\label{sec: S_continuous_formulation}
The steady Stokes problem in a two-dimensional parametrized domain $\Omega_o(\boldsymbol\mu)\subset\mathbb{R}^2$ read as: find $\boldsymbol{u}_o(\boldsymbol\mu): \Omega_o(\boldsymbol\mu) \to \mathbb{R}^2$ and $p_o(\boldsymbol\mu): \Omega_o(\boldsymbol\mu) \to \mathbb{R}^2$ such that
\begin{equation}
\begin{cases}
-\nu(\boldsymbol\mu)\Delta{\boldsymbol{u}_o}(\boldsymbol\mu)+\nabla{p}_o(\boldsymbol\mu)= \boldsymbol{0}& \text{ in } \Omega_o(\boldsymbol\mu),\  \\ 
\rm div\ {\boldsymbol{u}_o}(\boldsymbol\mu)=0 & \text{ in } \Omega_o(\boldsymbol\mu),\\
\boldsymbol{u}_o(\boldsymbol\mu)=\boldsymbol{g}_D(\boldsymbol\mu) & \text{ on } \Gamma_{D,o}(\boldsymbol\mu),\\
\boldsymbol{u}_o(\boldsymbol\mu)=\boldsymbol{0}, & \text{ on } \Gamma_{W,o}(\boldsymbol\mu),
\end{cases}
\label{eq:S_continuous}
\end{equation}
where $\boldsymbol{u}_o(\boldsymbol\mu)$ is the unknown velocity and $p_o(\boldsymbol\mu)$ is the unknown pressure, $\nu(\boldsymbol\mu)$ is the viscosity of fluid. These quantities, as well as the domain $\Omega_o(\boldsymbol\mu)$, may depend on a parameter $\boldsymbol\mu\in\mathbb{P}$ which accounts for parametrized physical properties or geometrical configurations. The boundary $\partial\Omega_o(\boldsymbol\mu)$ is divided into two parts in such a way that $\partial\Omega_o(\boldsymbol\mu)=\Gamma_{D,o}(\boldsymbol\mu)\cup\Gamma_{W,o}(\boldsymbol\mu),$ where $\Gamma_{D,o}(\boldsymbol\mu)$ is the Dirichlet boundary with non-homogeneous data and $\Gamma_{W,o}(\boldsymbol\mu)$ denotes the Dirichlet boundary with zero data (i.e., walls).
\par 
In order to write the weak formulation of problem (\ref{eq:S_continuous}), we introduce a reference domain, i.e. a $\boldsymbol\mu$-independent configuration $\Omega$ by assuming that each parametrized domain $\Omega_o(\boldsymbol\mu)$ can be obtained as the image of $\boldsymbol\mu$-independent domain $\Omega$ through a parametrized map $\boldsymbol{T}(.;\boldsymbol\mu):\mathbb{R}^2\rightarrow\mathbb{R}^2,$ i.e. $\Omega_o(\boldsymbol\mu)=\boldsymbol{T}(\Omega;\boldsymbol\mu).$
\par
Now the weak formulation of \eqref{eq:S_continuous} can be obtained by multiplying with the velocity, pressure test functions and using integration by parts; then by tracing everything back onto the reference domain $\Omega,$ we have the following parametrized weak formulation of problem \eqref{eq:S_continuous}:
\begin{equation}
\begin{cases}
\text{Find } \boldsymbol{u}(\boldsymbol\mu)\in \boldsymbol{V}, \, p(\boldsymbol\mu)\in Q:\\ a(\boldsymbol{u}(\boldsymbol\mu),\boldsymbol{v};\boldsymbol\mu)+b(\boldsymbol{v},p(\boldsymbol\mu);\boldsymbol\mu)=F(\boldsymbol{v};\boldsymbol\mu) & \forall \, v\in \boldsymbol{V},\\
b(\boldsymbol{u}(\boldsymbol\mu),q;\boldsymbol\mu)=G(q;\boldsymbol\mu) & \forall \, q\in Q,
\end{cases}
\label{eq:weak form1}
\end{equation}
where $\boldsymbol{V}=[H^1_0(\Omega)]^2$ and $Q=L^2_0(\Omega)=\lbrace q\in L^2(\Omega): \int_\Omega q=0\rbrace $.  Bilinear forms related to diffusion and pressure-divergence operators are defined as:
\begin{equation}
a(\boldsymbol{u},\boldsymbol{v};\boldsymbol\mu)=\int_{\Omega}\dfrac{\partial\boldsymbol{u}}{\partial x_i}\kappa_{ij}(x;\boldsymbol\mu)\dfrac{\partial\boldsymbol{v}}{\partial x_j}d\boldsymbol{x}, \qquad b(\boldsymbol{v},q;\boldsymbol\mu)=-\int_{\Omega}q \chi_{ij}(x;\boldsymbol\mu)\dfrac{\partial{v_j}}{\partial x_i}d\boldsymbol{x}.
\label{eq: S_bilinear_forms}
\end{equation}
The transformation tensors for bilinear viscous and pressure divergence terms in \eqref{eq: S_bilinear_forms} are defined as follows \cite{Ballarin2015}: 
\begin{equation}
\begin{split}
\boldsymbol{\kappa}(\boldsymbol{x};\boldsymbol\mu)&=\nu(\boldsymbol\mu)(J_T(\boldsymbol{x};\boldsymbol\mu))^{-1}(J_T(\boldsymbol{x};\boldsymbol\mu))^{-T}|J_T(\boldsymbol{x};\boldsymbol\mu)|,\\
\boldsymbol{\chi}(\boldsymbol{x};\boldsymbol\mu)&=(J_T(\boldsymbol{x};\boldsymbol\mu))^{-1}|J_T(\boldsymbol{x};\boldsymbol\mu)|,
\end{split}
\label{eq: S_tensors}
\end{equation}
where $|J_T|$ is the determinant of the Jacobian matrix $J_T\in \mathbb{R}^{2\times 2}$ of the map $T(.;\boldsymbol\mu)$. $F$ and $G$ are terms due to non-homogeneous Dirichlet boundary condition on the boundary are defined as:
\begin{equation}
\begin{split}
F(\boldsymbol{v};\boldsymbol\mu)&=-a(\boldsymbol{l(\mu)},\boldsymbol{v};\boldsymbol\mu),\\
G(q;\boldsymbol\mu)&=-b(\boldsymbol{l(\mu)},q;\boldsymbol\mu),
\end{split}
\label{eq: S_lift_functions}
\end{equation}
where we denote by $\boldsymbol{l}(\boldsymbol\mu)$ a parametrized lifting function such that  $\boldsymbol{l}(\boldsymbol\mu)|_{\Gamma_{D_g}}=\boldsymbol{g}_D(\boldsymbol\mu).$ 
\subsection{Finite Element formulation}
\label{subsec: FE_Stokes}
For a given parameter $\boldsymbol\mu\in \mathbb{P},$ the Galerkin-FE approximation of the parametrized Stokes problem \eqref{eq:weak form1} reads as follows:
\begin{equation}
\begin{cases}
\text{Find } \boldsymbol{u}_{h}(\boldsymbol{\mu})\in \boldsymbol{V}_{h}, p_{h}(\boldsymbol{\mu})\in Q_{h}:\\ 
a(\boldsymbol{u}_h(\boldsymbol\mu),\boldsymbol{v}_h;\boldsymbol\mu)+b(\boldsymbol{v}_h,p_h(\boldsymbol\mu);\boldsymbol\mu)=F(\boldsymbol{v}_h;\boldsymbol\mu) & \forall \, \boldsymbol{v}_h\in \boldsymbol{V}_h,\\ b(\boldsymbol{u}_h(\boldsymbol\mu),q_h;\boldsymbol\mu)=G(q_h;\boldsymbol\mu) & \forall \, q_h\in Q_h,
\end{cases}
\label{eq:weak discrete_stokes}
\end{equation}
where $\boldsymbol{V}_{h}$ and $Q_{h}$ are finite dimensional subspaces of $\boldsymbol{V}$ and $Q$ of dimension $\mathcal{N}_u$ and $\mathcal{N}_p$, respectively with $h$ related to the computational mesh size. When dealing with saddle point formulations it is well known that $\boldsymbol{V}_{h}$ and $Q_{h}$ should be judiciously chosen \cite{QV}; we will get back to this point at the end of this subsection.\par 
Let $\{\boldsymbol\phi_i^h\}_{i=1}^{\mathcal{N}_u}$ and $\{\psi_j^h\}_{j=1}^{\mathcal{N}_p}$ be basis functions of $\boldsymbol{V}_h$ and $Q_h$ respectively. We introduce the matrices $A(\boldsymbol\mu)\in\mathbb{R}^{\mathcal{N}_u\times\mathcal{N}_u}$ and $B(\boldsymbol\mu)\in\mathbb{R}^{\mathcal{N}_p\times\mathcal{N}_u}$ whose entries are
\begin{equation}
\left(A(\boldsymbol\mu)\right)_{ij}=a(\boldsymbol\phi_j^h,\boldsymbol\phi_i^h;\boldsymbol\mu), \quad \left(B(\boldsymbol\mu)\right)_{ki}=b(\boldsymbol\phi_i^h,\psi_k^h;\boldsymbol\mu), \text{ for } 1\leq i,j\leq\mathcal{N}_u, 1\leq k\leq\mathcal{N}_p,
\end{equation}
and the algebraic form of discrete problem \eqref{eq:weak discrete_stokes} problem reads
\begin{equation}
\left[\begin{array}{cc} A(\boldsymbol\mu) & B^T(\boldsymbol\mu)\\ B(\boldsymbol\mu) & \boldsymbol{0} \end{array}\right]
\left[\begin{array}{c} \boldsymbol{U}(\boldsymbol\mu) \\ \boldsymbol{P}(\boldsymbol\mu) \end{array}\right] 
= \left[\begin{array}{c} \boldsymbol{\bar{f}}(\boldsymbol\mu) \\ \boldsymbol{\bar{g}}(\boldsymbol\mu) \end{array}\right] 
\label{eq: Algebraic_Stokes}
\end{equation}
for the vectors $\boldsymbol{U}=(u_h^{(1)},...,u_h^{({\mathcal{N}_u})})^T, \boldsymbol{P}=(p_h^{(1)},...,p_h^{({\mathcal{N}_p})})^T$ such that $\boldsymbol{u}_h = \sum_{i=1}^{\mathcal{N}_u} u_h^{(i)} \boldsymbol\phi_i^h$ and $p_h = \sum_{j=1}^{\mathcal{N}_p} p_h^{(j)} \psi_j^h$, where for $1\leq i\leq\mathcal{N}_u$ and $1\leq k\leq\mathcal{N}_p$:
\begin{equation}
(\boldsymbol{\bar{f}}(\boldsymbol\mu))_i=-a(\boldsymbol{l}_h,\boldsymbol\phi_i^h;\boldsymbol\mu), \quad (\boldsymbol{\bar{g}}(\boldsymbol\mu))_k=-b(\boldsymbol{l}_h,\psi_k^h;\boldsymbol\mu), 
\end{equation}
with $\boldsymbol{l}_h=\boldsymbol{l}_h(\boldsymbol\mu)$, $\boldsymbol{l}(\boldsymbol\mu)$ in \eqref{sec: S_continuous_formulation} is discretized by the FE interpolant.
For the sake of the subsequent reduction procedure, affine parametric dependence is assumed for operators in \eqref{eq: Algebraic_Stokes} \cite{RB2016}:
\begin{equation}
\begin{split}
A(\boldsymbol\mu)&=\sum_{q=1}^{Q_a}\Theta_{q}^{a}(\boldsymbol\mu)A^q, \qquad B(\boldsymbol\mu)=\sum_{q=1}^{Q_b}\Theta_{q}^{b}(\boldsymbol\mu)B^{q}, \\
\boldsymbol{\bar{f}}(\boldsymbol\mu)&=\sum_{q=1}^{Q_f}\Theta_{q}^{f}(\boldsymbol\mu)\boldsymbol{\bar{f}}^q,\qquad  \boldsymbol{\bar{g}}(\boldsymbol\mu)=\sum_{q=1}^{Q_g}\Theta_{q}^{g}(\boldsymbol\mu)\boldsymbol{\bar{g}}^q.
\end{split}
\label{eq: affine_NS}
\end{equation}
\par 
Finally, coming back to the choice of the FE pairs, the FE spaces $\boldsymbol{V}_h$ and $Q_h$  have to fulfill the following parametrized version of the inf-sup condition \cite{QV}:
\begin{equation} 
\exists \beta_{0}(\boldsymbol{\mu})>0: \beta_{h}(\boldsymbol{\mu})= \inf_{q_h\in Q_{h}}\sup_{\boldsymbol{v}_h\in \boldsymbol{V}_{h}}\frac{b(\boldsymbol{v}_h,q_h;\boldsymbol\mu)}{\|\boldsymbol{v}_h\|_{\boldsymbol{V}_h}\|q_h\|_{Q_h}}\geq \beta_{0}(\boldsymbol{\mu}) \quad\forall \boldsymbol\mu\in \mathbb{P}.
\label{eq: inf-sup}
\end{equation}
This relation holds if, e.g., the Taylor-Hood $(\mathbb{P}_2/\mathbb{P}_1)$ FE spaces are chosen. Condition \eqref{eq: inf-sup} does not hold in case of equal order FE spaces  $(\mathbb{P}_k/\mathbb{P}_k), k\geq 1$ and for lowest order element $(\mathbb{P}_1/\mathbb{P}_0)$. To handle these cases, a stabilized formulation is introduced as in the following.\nocite{Rozza2009}
\subsection{Stabilized Finite Element formulation}
\label{subsec: FE_Stokes_stab}
When the finite dimensional spaces $\boldsymbol{V}_{h}$ and $Q_{h}$ do not satisfy the $\inf$-$\sup$ stability condition (\ref{eq: inf-sup}), then in order to avoid possible spurious pressure modes the use of a stabilized formulation, introduced in \cite{Hughes1986, Brooks1982, Hughes1989, Hughes1980, Douglas1989, ZHU1993}, is nowadays widespread in the FE community. This amounts to modifying \eqref{eq:weak discrete_stokes} as 
\begin{equation}
\begin{cases}
\text{Find } \boldsymbol{u}_{h}(\boldsymbol{\mu})\in \boldsymbol{V}_h, \, p_{h}(\boldsymbol{\mu})\in Q_{h}:\\ a(\boldsymbol{u}_{h}(\boldsymbol{\mu}),\boldsymbol{v}_{h};\boldsymbol{\mu})+b(\boldsymbol{v}_{h},p_{h}(\boldsymbol{\mu});\boldsymbol{\mu})-s^{u,v}_{h}(\boldsymbol{u}_{h},\boldsymbol{v}_{h};\boldsymbol\mu)-s^{p,v}_{h}(p_h,\boldsymbol{v}_{h};\boldsymbol\mu)=F(\boldsymbol{v}_h;\boldsymbol\mu) & \forall \, \boldsymbol{v}_{h}\in \boldsymbol{V}_{h},\\ b(\boldsymbol{u}_{h}(\boldsymbol{\mu}),q_{h};\boldsymbol{\mu})-s^{u,q}_{h}(\boldsymbol{u}_{h},q_{h};\boldsymbol\mu)-s^{p,q}_{h}(p_h,q_{h};\boldsymbol\mu)=G(q_h;\boldsymbol\mu) & \forall \, q_{h}\in Q_{h},
\end{cases}
\label{eq:stable_discrete}
\end{equation}
where $s^{u,v}_{h}(.,.;\boldsymbol\mu)$, $s^{p,v}_{h}(.,.;\boldsymbol\mu)$, $s^{u,q}_{h}(.,.;\boldsymbol\mu)$ and $s^{p,q}_{h}(.,.;\boldsymbol\mu)$ are the stabilization terms \cite{QV}. 
Following a strongly consistent residual based approach, these stabilization terms are defined as
\begin{equation}
s^{u,v}_{h}(\boldsymbol{u}_{h},\boldsymbol{v}_{h};\boldsymbol\mu):=\delta\sum_{K}h_K^{2}\int_K(-\nu\Delta\boldsymbol{u}_{h},-\rho\nu\Delta\boldsymbol{v}_{h}),
\label{eq:choice_general}
\end{equation}
\begin{equation}
s^{p,v}_{h}(p_h,\boldsymbol{v}_{h};\boldsymbol\mu):=\delta\sum_{K}h_K^{2}\int_K(\nabla{p}_{h},-\rho\nu\Delta\boldsymbol{v}_{h}),
\label{eq:choice_general_1}
\end{equation}
\begin{equation}
s^{u,q}_{h}(\boldsymbol{u}_{h},q_{h};\boldsymbol\mu):=\delta\sum_{K}h_K^{2}\int_K(-\nu\Delta\boldsymbol{u}_{h},\nabla{q}_{h}),
\label{eq:our_choice}
\end{equation}
and
\begin{equation}
s^{p,q}_{h}(p_h,q_{h};\boldsymbol\mu):=\delta\sum_{K}h_K^{2}\int_K(\nabla{p}_{h},\nabla{q}_{h}),
\label{eq:our_choice_1}
\end{equation}
where $K$ is an element of the triangulation of the reference domain $\Omega$, $h_K$ is the diameter of element $K$, $\delta$ is the stabilization coefficient assumed to be constant. For $\rho=0,1,-1$, the method (\ref{eq:choice_general})-\eqref{eq:our_choice_1} is respectively known as the pressure-poisson stabilized Galerkin (Franca-Hughes) \cite{Hughes1986}, Galerkin least-squares (GALS) \cite{HUGHES1987}, Douglas-Wang (DW) \cite{Douglas1989}. In case of linear interpolation for velocity and pressure $(\mathbb{P}_{1}/\mathbb{P}_{1})$ the Laplacian term $-\nu\Delta\boldsymbol{u}_{h}$ inside the stabilization vanishes and all above choices reduce to Brezzi-Pitk\"{a}ranta stabilization \cite{Brezzi1984}, written as
\begin{equation}
s^{u, v}_{h}(\boldsymbol{u}_{h},\boldsymbol{v}_{h};\boldsymbol\mu)=s^{p, v}_{h}(p_{h},\boldsymbol{v}_{h};\boldsymbol\mu)=s^{u,q}_{h}(\boldsymbol{u}_{h},q_h;\boldsymbol\mu):= 0, \qquad s^{p,q}_{h}(p_h,q_{h};\boldsymbol\mu):=\delta\sum_{K}h_K^{2}\int_K\nabla{p}_{h}\cdot\nabla{q}_{h}.
\label{eq:choice_BP}
\end{equation}
Further possible stabilization options are available, see e.g. \cite[chapter 9]{QV} and references therein.
Dependence of these stabilization terms on parameter $\boldsymbol\mu$ is motivated by the dependence through the parametrized solution. Further dependence on $\boldsymbol\mu$ might be warranted for the stabilization coefficient $\delta$ (see e.g. numerical test cases in \cite{Torlo2018} for the case of geometrical parametrization). In particular, it is often proposed in the FE literature (see e.g. \cite{Jenkins2013}) that the stabilization coefficient depends nonlinearly on the solution. For the sake of guaranteeing affine dependence on the stabilization terms, here we assume instead $\delta$ to be a constant, with no noticeable deterioration in the numerical results provided that the stabilization coefficient $\delta$ has been properly chosen. It cannot be too small otherwise the stabilization will be poor and spurious modes will not be eliminated, while a large value of parameter $\delta$ could result in a poor approximation for the pressure field near to the boundary.
Beyond this, if needed, parameter and solution dependent stabilization coefficient could be dealt with the empirical interpolation method \cite{Barrault2004}.
\par 
For the sake of exposition, in this work we will only focus on the stabilization option corresponding to $\rho=0$ (Franca-Hughes) \cite{Hughes1986}, i.e, the stabilization terms \eqref{eq:choice_general} and \eqref{eq:choice_general_1} are not taken into account.
\par 
After adding the stabilization terms into the system \eqref{eq: Algebraic_Stokes}, the stabilized algebraic formulation reads
\begin{equation}
\left[\begin{array}{cc} A(\boldsymbol\mu) & B^T(\boldsymbol\mu)\\ \tilde{B}(\boldsymbol\mu) & -S(\boldsymbol\mu) \end{array}\right]
\left[\begin{array}{c} \boldsymbol{U}(\boldsymbol\mu) \\ \boldsymbol{P}(\boldsymbol\mu) \end{array}\right] 
= \left[\begin{array}{c} \boldsymbol{\bar{f}}(\boldsymbol\mu) \\ \boldsymbol{\bar{g}}(\boldsymbol\mu) \end{array}\right] 
\label{eq: Algebraic_Stokes_stab},
\end{equation}
where $\tilde{B}(\boldsymbol\mu)$ and $S(\boldsymbol\mu)$ contains the effects of stabilization and defined as follows:
\begin{equation}
\left(\tilde{B}(\boldsymbol\mu)\right)_{ki}=b(\boldsymbol\phi_i^h,\psi_k^h;\boldsymbol\mu)+s^{u,q}_{h}(\boldsymbol\phi_i^h,\psi_k^h;\boldsymbol\mu),\quad
(S(\boldsymbol\mu))_{ij}=s^{p,q}_{h}(\psi^{h}_j,\psi^{h}_i;\boldsymbol\mu), \text{ for } 1\leq i,j\leq\mathcal{N}_u, 1\leq k\leq\mathcal{N}_p,
\end{equation}
\par
The motivation to choose a \textit{strongly consistent residual based method} is, it improve the stability of Galerkin FE method without compromising the consistency. This choice also allows the stabilized formulation to fulfill a modified version of \eqref{eq: inf-sup} with a supplementary terms \cite{boffi2013mixed,Burman2009,Becker2001}. Indeed, e.g. for the Brezzi-Pitk\"{a}ranta method, the stabilized formulation requires the FE spaces to fulfill the following modified $\inf$-$\sup$ condition
\begin{equation}
\exists \beta_{0}(\boldsymbol{\mu})>0:  \sup_{{\boldsymbol{v}_h}\in \boldsymbol{V}_{h}}\frac{b(\boldsymbol{v}_h,q_h;\boldsymbol\mu)}{\|\nabla{\boldsymbol{v}_h}\|}+s^{p,q}_{h}\left(q_{h},q_{h};\boldsymbol\mu\right)^{1/2}\geq \beta_0(\boldsymbol\mu)\|q_h\|, \forall q_h\in Q_{h},
\label{eq:modified_inf-sup}
\end{equation}
This condition is a crucial motivation of our work, and will be adapted to the reduced order system in subsection \ref{subsec:RBStokes_stab}.
\subsection{Reduced Basis formulation}
\label{subsec:RBStokes}
In this section we present the RB formulation of steady Stokes problem. The RB method seeks the approximation of FE solution to \eqref{eq:weak discrete_stokes}. In first step we construct a set of global basis functions. Let $\boldsymbol\mu^1,...,\boldsymbol\mu^N$ be a set of snapshot parameter values chosen by greedy algorithm \cite{Rozza2013}. We denote by $\boldsymbol{u}_h(\boldsymbol{\mu}^n)$ and $p_h(\boldsymbol{\mu}^n)$, $n = 1, \hdots, N$, the corresponding snapshot solutions
for the velocity and pressure. The snapshot solutions are obtained as solutions of a stable FE formulation of \eqref{eq:weak discrete_stokes}. As our focus (especially in the nonlinear case in the second part of the work) is not on the certification of the reduced model, during the greedy iterations we rely on an error indicator based on the residual, rather than a certified error bound \cite{RB2016,Veroy2012} which would need to be properly extended to the stabilized case. Nonetheless, the availability of error bounds does not affect the presentation of the rest of the methodology described in this paper.
At the end of the greedy procedure, we obtain the reduced velocity space $\boldsymbol{V}_N\subset \boldsymbol{V}_{h}$ and reduced pressure space $Q_N\subset Q_{h}$, respectively as:
\begin{equation}
\boldsymbol{V}_N = \text{span}\left\lbrace\boldsymbol{u}_{h}(\boldsymbol\mu^n), 1\leq n\leq N_u\right\rbrace,
\end{equation}
and 
\begin{equation}
Q_N = \text{span}\left\lbrace p_{h}(\boldsymbol\mu^n), 1\leq n\leq N_p\right\rbrace, 
\end{equation}
where $N_u = N_p = N$ are the dimensions of RB velocity space $\boldsymbol{V}_N$ and RB pressure space $Q_N$, respectively. Applying the Gram-Schmidt orthogonalization process on the snapshots \cite{RB2016}, we denote by $\{\boldsymbol\xi^u_n\}^{N_u}_{n=1}$ and $\{\xi^p_n\}^{N_p}_{n=1}$ mutually orthonormal functions obtained from velocity and pressure snapshots, respectively, which we are going to use as basis functions for the reduced spaces instead of the snapshots.\par
As well known in the RB community, it is important to point out that, even when the spaces $\boldsymbol{V}_N$ and $Q_N$ are obtained collecting snapshots from a stable full order model from section 2.1, a Galerkin projection over the reduced spaces does not guarantee the fulfillment of the following reduced inf-sup condition:
\begin{equation}
\exists \beta_{0,N}(\boldsymbol{\mu})>0: \quad  \beta_N(\boldsymbol{\mu})= \inf_{q_N\in Q_N} \sup_{\boldsymbol{v}_N\in \boldsymbol{V}_N}\frac{b(\boldsymbol{v}_N,q_N;\boldsymbol{\mu})}{\|\boldsymbol{v}_N\|_{\boldsymbol{V}_N}\|q_N\|_{Q_N}}\geq \beta_{0,N}(\boldsymbol{\mu}) \quad \forall \boldsymbol\mu\in \mathbb{P}.
\label{eq: inf-sub RB}
\end{equation}
Indeed, following \cite{Veroy2007} the usual reduced order methodology relies on the enrichment of the RB velocity space with supremizer solutions.
Even though several variants have been discussed in literature \cite{Veroy2007, Rovas2003, Rozza2005a}, for the sake of exposition we will focus only on the so-called \emph{approximate supremizer enrichment}. This approach requires the introduction of the supremizer operator $T^{\boldsymbol\mu}:Q_h\rightarrow\boldsymbol{V}_h$ defined as follows:
\begin{equation}
(T^{\boldsymbol\mu}q_h,\boldsymbol{v}_h)_{\boldsymbol{V}}=b(\boldsymbol{v}_h,q_h;\boldsymbol\mu), \quad \forall \boldsymbol{v}\in \boldsymbol{V}_h.
\label{eq: supremizer}
\end{equation}
which is evaluated for $\boldsymbol\mu = \boldsymbol\mu^n$ and the corresponding pressure snapshot $q_h := p_{h}(\boldsymbol\mu^n)$, $n = 1, \hdots, N$, to obtain $N$ supremizer snapshots.
Afterwards, the RB velocity space $\boldsymbol{V}_N$ is enriched with the supremizer snapshots. We denote the enriched RB velocity space by $\tilde{\boldsymbol{V}}_N$, defined as:
\begin{equation}
\tilde{\boldsymbol{V}}_N = \text{span}\left\lbrace\boldsymbol{u}_{h}(\boldsymbol\mu^n), 1\leq n\leq N_u;T^{\boldsymbol\mu^n}p_{h}(\boldsymbol\mu^n), 1\leq n\leq N_s\right\rbrace, 
\end{equation}
where $N_s \leq N_p$ denotes the number of supremizer snapshots. The dimension of $\tilde{\boldsymbol{V}}_N$ will thus be $N_u+N_s$; in order to keep the notations simple, in the following we will always take $N_u=N_s=N_p=N$.
Now the RB formulation corresponding to FE problem \eqref{eq:weak discrete_stokes} can be written as:
\begin{equation}
\begin{cases}
\text{Find } \boldsymbol{u}_N(\boldsymbol{\mu}) \in \tilde{\boldsymbol{V}}_N, p_N(\boldsymbol{\mu})) \in Q_N:\\ a(\boldsymbol{u}_N(\boldsymbol{\mu}),\boldsymbol{v}_N;\boldsymbol{\mu})+b(\boldsymbol{v}_N,p_N(\boldsymbol{\mu}))=F(\boldsymbol{v}_N;\boldsymbol\mu) & \forall \, \boldsymbol{v}_N\in \tilde{\boldsymbol{V}}_N, \\ b(\boldsymbol{u}_N(\boldsymbol{\mu}),q_N;\boldsymbol{\mu})= G(q_N;\boldsymbol\mu) & \forall \, q_N\in Q_N.
\end{cases}
\label{eq:RBFormulation}
\end{equation}
The solution $(\boldsymbol{u}_N(\boldsymbol{\mu}),p_N(\boldsymbol{\mu}))\in \tilde{\boldsymbol{V}}_N\times Q_N$ of \eqref{eq:RBFormulation} can be expressed as a linear combination of the basis functions:
\begin{equation}
\boldsymbol{u}_N(\boldsymbol{\mu})=\sum_{n=1}^{2N}U_{N,n}(\boldsymbol{\mu})\boldsymbol\xi^u_n, \qquad p_N(\boldsymbol{\mu})=\sum_{n=1}^{N}P_{N,n}(\boldsymbol{\mu})\xi^p_n,
\label{eq:RBsolution}
\end{equation}
where $\boldsymbol{U}_N(\boldsymbol\mu) = [U_{N,n}(\boldsymbol{\mu})]_{n=1}^{N_u}$ and $\boldsymbol{P}_N(\boldsymbol\mu) = [P_{N,n}(\boldsymbol{\mu})]_{n=1}^{N_p}$ denote the vector of coefficients of the reduced basis approximation for velocity and pressure.
Finally, we write the system in compact form as 
\begin{equation}
\left[\begin{array}{cc} A_{N}(\boldsymbol\mu) & B^T_{N}(\boldsymbol\mu)\\ B_{N}(\boldsymbol\mu) & \boldsymbol{0} \end{array}\right]
\left[\begin{array}{c} \boldsymbol{U}_N(\boldsymbol\mu) \\ \boldsymbol{P}_N(\boldsymbol\mu) \end{array}\right] 
= \left[\begin{array}{c} \boldsymbol{\bar{f}}_N(\boldsymbol\mu) \\ \boldsymbol{\bar{g}}_N(\boldsymbol\mu) \end{array}\right], 
\label{eq: Algebraic_Stokes_RB}
\end{equation}
where the RB tensors are computed as 
\begin{equation}
\begin{split}
A_N(\boldsymbol\mu)&=Z^T_{u,s}A(\boldsymbol\mu)Z_{u,s}, \quad B_N(\boldsymbol\mu)=Z^T_{p}B(\boldsymbol\mu)Z_{u,s}, \\ \boldsymbol{\bar{f}}_N(\boldsymbol\mu)&=Z^T_{u,s}\boldsymbol{\bar{f}}(\boldsymbol\mu),\qquad
\boldsymbol{\bar{g}}_N(\boldsymbol\mu)=Z^T_{p}\boldsymbol{\bar{g}}(\boldsymbol\mu),
\end{split}
\label{eq:RBmatrices_SteadyStokes}
\end{equation}
being $Z_{u,s} \in \mathbb{R}^{\mathcal{N}_u \times N_{u,s}}$ and $Z_p \in \mathbb{R}^{\mathcal{N}_p \times N_{p}}$ rectangular matrices that contain the FE degrees of freedom of the basis of $\tilde{\boldsymbol{V}}_N$ and $Q_N$, respectively.
\subsection{Stabilized Reduced Basis formulation}
\label{subsec:RBStokes_stab}
We introduce now the stabilized Reduced Basis model derived from the stabilized FE problem \eqref{eq:stable_discrete}. During the offline phase, the computation of the reduced spaces is done in section 2.3, resulting in a reduced velocity space $\boldsymbol{V}_N$, enriched reduced velocity space $\tilde{\boldsymbol{V}}_N$, and reduced pressure space $Q_N$. During the online stage, the stabilized RB problem reads
\begin{equation}
\begin{cases}
\text{Find } \boldsymbol{u}_{N}(\boldsymbol{\mu})\in \boldsymbol{V}_N, \, p_{N}(\boldsymbol{\mu})\in Q_{N}:\\ a(\boldsymbol{u}_{N}(\boldsymbol{\mu}),\boldsymbol{v}_{N};\boldsymbol{\mu})+b(\boldsymbol{v}_{N},p_{N}(\boldsymbol{\mu});\boldsymbol{\mu})-s^{u,v}_{N}(\boldsymbol{u}_{N},\boldsymbol{v}_{N};\boldsymbol\mu)-s^{p,v}_{N}(p_N,\boldsymbol{v}_{N};\boldsymbol\mu)=F(\boldsymbol{v}_N;\boldsymbol\mu) & \forall \, \boldsymbol{v}_{N}\in \boldsymbol{V}_{N},\\ b(\boldsymbol{u}_{N}(\boldsymbol{\mu}),q_{N};\boldsymbol{\mu})-s^{u,q}_{N}(\boldsymbol{u}_{N},q_{N};\boldsymbol\mu)-s^{p,q}_{N}(p_N,q_{N};\boldsymbol\mu)=G(q_N;\boldsymbol\mu) & \forall \, q_{N}\in Q_{N},
\end{cases}
\label{eq:RBFormulation_stab}
\end{equation}
where $s^{u,v}_{N}(.,.;\boldsymbol\mu)$, $s^{p,v}_{N}(.,.;\boldsymbol\mu)$, $s^{u,q}_{N}(.,.;\boldsymbol\mu)$ and $s^{p,q}_{N}(.,.;\boldsymbol\mu)$ are the reduced order stabilization terms defined as
\begin{equation}
s^{u,v}_{N}(\boldsymbol{u}_{N},\boldsymbol{v}_{N};\boldsymbol\mu):=\delta\sum_{K}h_K^{2}\int_K(-\nu\Delta\boldsymbol{u}_{N},-\rho\nu\Delta\boldsymbol{v}_{N}),
\label{eq:choice_general_RB}
\end{equation}
\begin{equation}
s^{p,v}_{N}(p_h,\boldsymbol{v}_{N};\boldsymbol\mu):=\delta\sum_{K}h_K^{2}\int_K(\nabla{p}_{N},-\rho\nu\Delta\boldsymbol{v}_{N}),
\label{eq:choice_general_1_RB}
\end{equation}
\begin{equation}
s^{u,q}_{N}(\boldsymbol{u}_{N},q_{N};\boldsymbol\mu):=\delta\sum_{K}h_K^{2}\int_K(-\nu\Delta\boldsymbol{u}_{N},\nabla{q}_{N}),
\label{eq:our_choice_RB}
\end{equation}
and
\begin{equation}
s^{p,q}_{N}(p_N,q_N;\boldsymbol\mu):=\delta\sum_{K}h_K^{2}\int_K(\nabla{p}_{N},\nabla{q}_{N}),
\label{eq:our_choice_1_RB}
\end{equation}
In a similar way, one could seek solutions in the enriched velocity space by replacing $\boldsymbol{V}_N$ with $\tilde{\boldsymbol{V}}_N$ in \eqref{eq:RBFormulation_stab}. With slight abuse of notation, we will keep denoting by $(\boldsymbol{u}_N(\boldsymbol{\mu}),p_N(\boldsymbol{\mu}))$ solutions to either \eqref{eq:RBFormulation} and \eqref{eq:RBFormulation_stab}, as it will be clear from the context to which RB formulation we will refer to. The corresponding algebraic formulation reads
\begin{equation}
\left[\begin{array}{cc} A_N(\boldsymbol\mu) & B_N^T(\boldsymbol\mu)\\ \tilde{B}_N(\boldsymbol\mu) & -S_N(\boldsymbol\mu) \end{array}\right]
\left[\begin{array}{c} \boldsymbol{U}_N(\boldsymbol\mu) \\ \boldsymbol{P}_N(\boldsymbol\mu) \end{array}\right] 
= \left[\begin{array}{c} \boldsymbol{\bar{f}}_N(\boldsymbol\mu) \\ \boldsymbol{\bar{g}}_N(\boldsymbol\mu) \end{array}\right] 
\label{eq: Algebraic_Stokes_stab_RB},
\end{equation}
where $\tilde{B}_N(\boldsymbol\mu)$ and $S_N(\boldsymbol\mu)$ are RB stabilization matrices defined as:
\begin{equation}
\tilde{B}_N(\boldsymbol\mu)=Z^T_{p}\tilde{B}(\boldsymbol\mu)Z_{u,s}, \quad S_N(\boldsymbol\mu)=Z^T_{p}S(\boldsymbol\mu)Z_{p}, 
\label{eq:RBmatrices_SteadyStokes_STAB}
\end{equation}
To further present our methodology, we also define the reduced order version of modified $\inf$-$\sup$ condition \eqref{eq:modified_inf-sup} as:
\begin{equation}
\exists \beta_{0,N}(\boldsymbol{\mu})>0: \sup_{{\boldsymbol{v}_N}\in \boldsymbol{V}_{N}}\frac{b(\boldsymbol{v}_N,q_N;\boldsymbol\mu)}{\|\nabla{\boldsymbol{v}_N}\|}+s^{p,q}_{N}(q_{N},q_{N};\boldsymbol\mu)^{1/2}\geq \beta_{0,N}(\boldsymbol{\mu}) \|q_N\|, \forall q_N\in Q_{N},
\label{eq:modified_inf-sup_RB}
\end{equation}
where $s^{p,q}_{N}(.,.;\boldsymbol\mu)$ is due to the addition of stabilization terms in RB formulation. It is now clear that the two addend on the left-hand side can contribute in different ways to the overall inf-sup stability of the reduced problem. In particular, one could either increase the first term by replacing $\boldsymbol{V}_{N}$ with $\tilde{\boldsymbol{V}}_N$, thus exploiting the existing supremizer enrichment procedure, or rely on the contribution of the second term due to the underlying FE residual based stabilization, or both. Overall, this results in the following four different combinations:
\begin{enumerate}
\item[\textit{(i)}] the first option is to increase both the first addend, through RB velocity space enrichment with supremizer solutions, and the second addend as well, through residual based stabilization. Thus, during the online stage we will solve \eqref{eq:RBFormulation_stab}, upon replacing $\boldsymbol{V}_N$ with $\tilde{\boldsymbol{V}}_N$. As the residual based stabilization terms are present both \textit{offline} and \textit{online}, we will call this methodology \textit{offline-online stabilization} in agreement with the denomination introduced in \cite{PR014a}. Furthermore, as supremizer enrichment is performed, we will denote this option by \textit{offline-online stabilization} with supremizer;
\item[\textit{(ii)}] as a second option, we do not enrich the RB velocity space with supremizer solutions, yet we still rely on the residual based stabilization during the online phase. Thus, during the online stage, we solve \eqref{eq:RBFormulation_stab}. When compared to (i), this option seems attractive in terms of CPU time because it decreases the dimensionality of the velocity space, resulting in possible larger speedups. We call this option \textit{offline-online stabilization} without supremizer; 
\item[\textit{(iii)}] the third option is to enrich the RB velocity space with supremizer solutions and neglect stabilization during the online phase. Thus, online we end up solving \eqref{eq:RBFormulation}. If compared to (i), this option seems attractive in perspective for cases where the stabilization coefficient $\delta$ depends nonlinearly on the solution, and an accurate approximation of the stabilization terms would require preprocessing through EIM. We denote this case with \textit{offline-only stabilization} with supremizer, the first part of the name underlining the fact that stabilization is applied only offline;
\item[\textit{(iv)}] a fourth (fictitious) option could be to completely avoid both supremizer enrichment and residual based stabilization. We would call this option as \textit{offline-only stabilization} without supremizer. It is clear, however, that this methodology would hardly grant a positive reduced inf-sup constant, as neither stabilizing contributions are enabled. For this reason, results for this case are not reported in numerical results presented in this paper.
\end{enumerate}
In the following, we employ equal order FE spaces  $(\mathbb{P}_k/\mathbb{P}_k), k\geq 1$ during the offline stage. We remark that, during this stage, reduced solutions required by the greedy algorithm are computed via option (i). Indeed, the consistency of the \textit{offline-online stabilization} is necessary to ensure that the same parameter will not be selected twice by the greedy algorithm. In contrast, during the online stage we study the performance of options (i)-(iii).
\section{Numerical results and discussion for parametrized Stokes problems}
\label{sec: Numerical_Results}
In this section, we present some numerical results for stabilized reduced order model for steady Stokes problem developed in section \ref{sec: S_continuous_formulation}. Numerical simulations are carried out in FreeFem++ \cite{Hetch2013}, and also with RBniCS \cite{Ballarin2016rb} for comparison.
\par
\begin{figure}[H]
\centering
\includegraphics[width=0.45\textwidth]{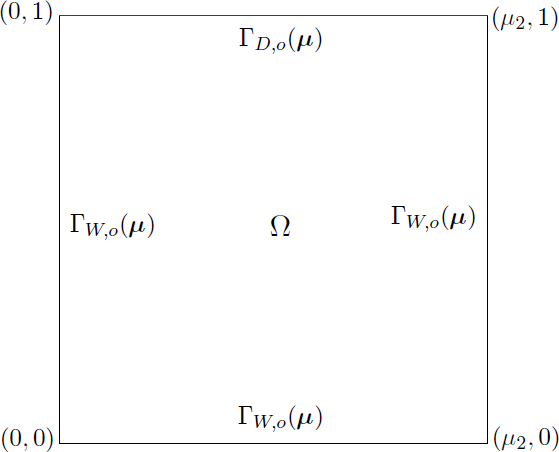}
\caption{Parametrized domain.}
\label{fig:Domain}
\end{figure}

 
As a test case we consider the parametrized cavity flow problem \cite{Hughes1986}. We set the parametrized domain $\Omega_o(\boldsymbol\mu)=(0,1+\mu_2)\times(0,1)$, where we define $\boldsymbol{\mu}=(\mu_1,\mu_2)$ such that $\mu_1$ is a physical parameter (kinematic viscosity of fluid) and $\mu_2$ is a geometrical parameter (length of domain). The parametrized domain is shown in Fig. \ref{fig:Domain}, along with the partition $\Gamma_{D,o}(\boldsymbol\mu) \cup \Gamma_{W,o}(\boldsymbol\mu)$ of its boundary $\partial\Omega$; unit horizontal velocity is imposed on the side $\Gamma_{D,o}(\boldsymbol\mu)$.
\subsection{Numerical results for $\mathbb{P}_k/\mathbb{P}_k$ for $k=1,2$}
\label{subsec:StokesP1P1}
Since a stabilized FE formulation is used, we need not to employ Taylor-Hood elements, rather we can rely on equal order velocity-pressure FE pairs. The results in this section are obtained for $\mathbb{P}_k/\mathbb{P}_k$ for $k=1,2$. The case $k=1$ is attractive in the FE context as it allows to employ the lowest order continuous FE spaces for both velocity and pressure. However, unless suitable gradient recovery techniques are used, all the various residual based stabilization techniques \eqref{eq:choice_general}-\eqref{eq:our_choice} coincide with the simplest Brezzi-Pitk{\"a}ranta stabilization \eqref{eq:choice_BP} when using lowest order elements. Therefore, in order to test the applicability of the Franca-Hughes stabilization we also provide results for $k=2$. Details of the problem and computational costs are summarized in Table \ref{tab:comput_steady_stokes_GEO}.\par
\begin{figure}[H]
  \centering
  \begin{subfigure}{0.48\textwidth}
    \centering\includegraphics[width=\textwidth]{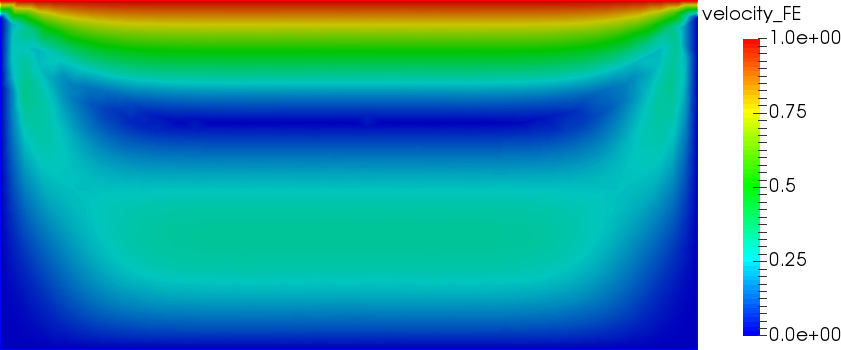}
    \caption{FE Velocity}
  \end{subfigure}%
  \quad%
  \begin{subfigure}{0.48\textwidth}
    \centering\includegraphics[width=\textwidth]{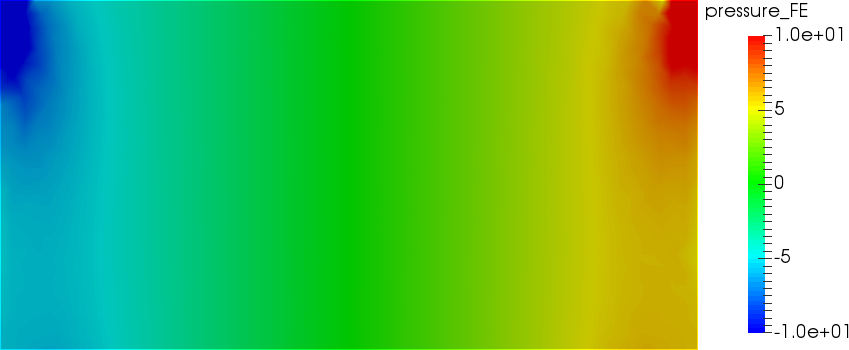}
    \caption{FE Pressure}
  \end{subfigure}
    \quad%
   \begin{subfigure}{0.48\textwidth}
    \centering\includegraphics[width=\textwidth]{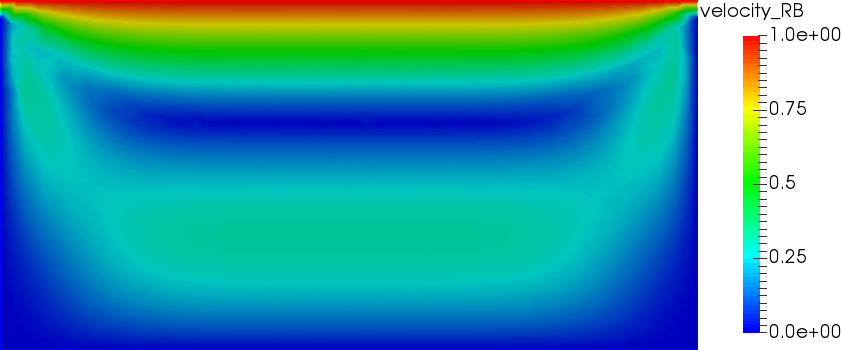}
    \caption{RB Velocity: \textit{offline-online stabilization} with supremizer}
  \end{subfigure}%
  \quad%
  \begin{subfigure}{0.48\textwidth}
    \centering\includegraphics[width=\textwidth]{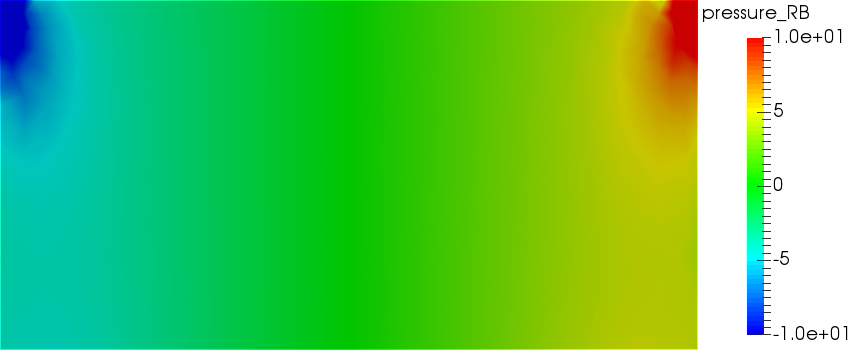}
    \caption{RB Pressure: \textit{offline-online stabilization} with supremizer}
  \end{subfigure}
    \quad%
    \begin{subfigure}{0.48\textwidth}
    \centering\includegraphics[width=\textwidth]{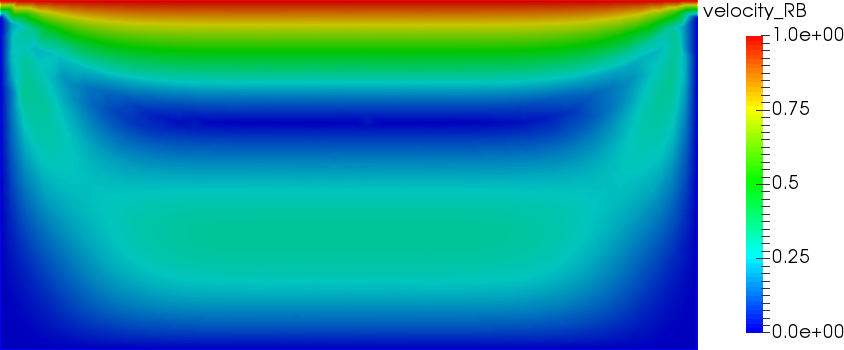}
    \caption{RB Velocity: \textit{offline-online stabilization} without supremizer}
  \end{subfigure}%
  \quad%
  \begin{subfigure}{0.48\textwidth}
    \centering\includegraphics[width=\textwidth]{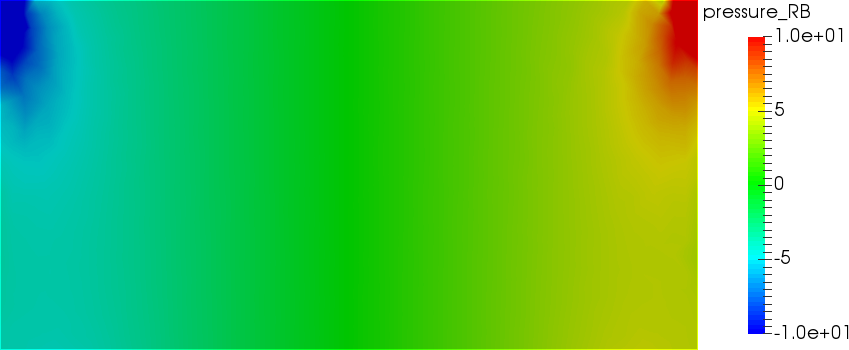}
    \caption{RB Pressure: \textit{offline-online stabilization} without supremizer}
  \end{subfigure}
  \quad%
   \begin{subfigure}{0.48\textwidth}
    \centering\includegraphics[width=\textwidth]{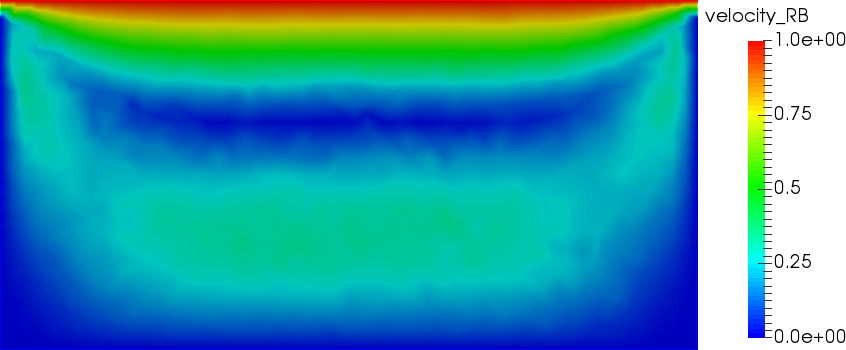}
    \caption{RB Velocity: \textit{offline-only} supremizer}
  \end{subfigure}%
  \quad%
  \begin{subfigure}{0.48\textwidth}
    \centering\includegraphics[width=\textwidth]{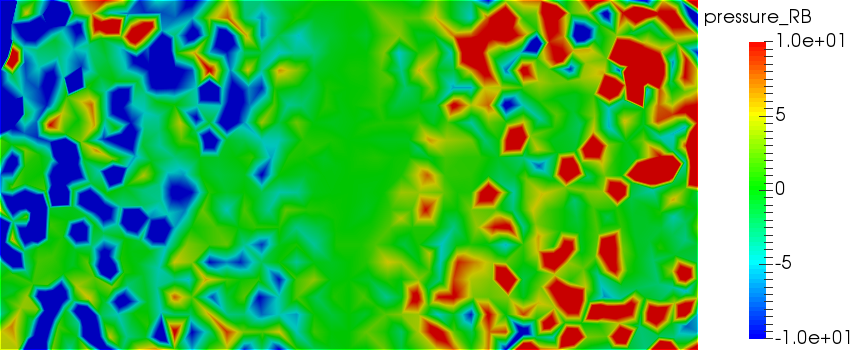}
    \caption{RB Pressure: \textit{offline-only} supremizer}
  \end{subfigure}
    \quad%
  \caption{Stokes problem: FE and RB solutions for velocity and pressure at $(\mu_1,\mu_2)=(0.6,2)$; $N_u=N_p=20$ and using $\mathbb{P}_2/\mathbb{P}_2$.}
\label{fig:stokes Sol}
\end{figure}

In Fig. \ref{fig:stokes Sol} we show a comparison between the FE velocity solution and the RB solutions obtained for three different options (i)-(iii), for $(\mu_1,\mu_2)=(0.6,2)$. These plots are shown in the $\mathbb{P}_2/\mathbb{P}_2$ case, and the corresponding plots for the $\mathbb{P}_1/\mathbb{P}_1$ pair hardly look any different.
From Fig. \ref{fig:stokes Sol}, we see that the RB velocity and pressure solutions obtained by using the \textit{offline-online stabilization} with/without supremizer looks similar to the FE solution. However, the RB solutions obtained by the \textit{offline-only stabilization} is poor, in particular pressure solution is highly oscillatory.\par
In order to see a more quantitive comparison between the \textit{offline-online stabilization} with/without supremizer and \textit{offline-only stabilization} with supremizer we report next the results of an error analysis, starting from the $\mathbb{P}_1/\mathbb{P}_1$ FE pair with $\delta = 0.05$.
These results show that \textit{offline-online stabilization} result in the most accurate methods, while \textit{offline-only stabilization} is inaccurate. We also see that the enrichment of supremizer together with \textit{offline-online stabilization} is beneficial for the pressure approximation, improving results by more than an order of magnitude, while has negigible effect on the velocity. \par
\begin{figure}[H]
  \centering
  \begin{subfigure}{0.48\textwidth}
    \centering\includegraphics[width=\textwidth]{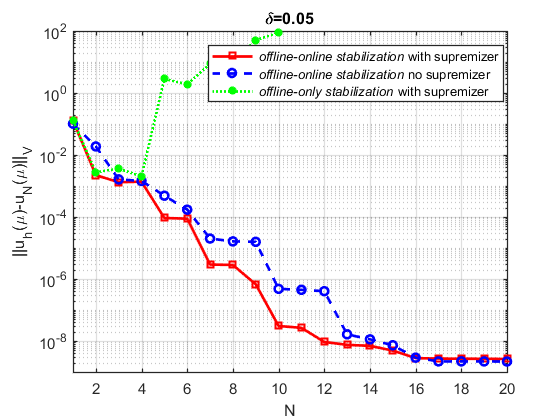}   
  \end{subfigure}%
  \quad%
  \begin{subfigure}{0.48\textwidth}
    \centering\includegraphics[width=\textwidth]{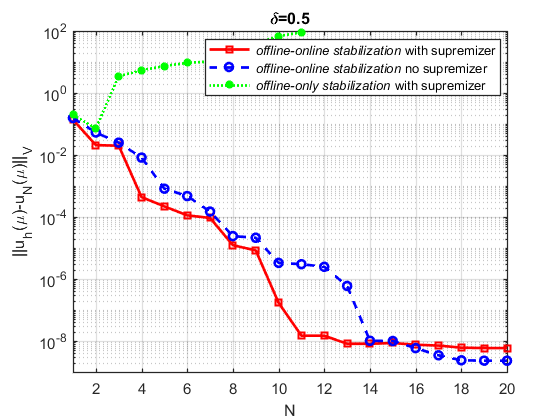}
  \end{subfigure}
  \caption{Stokes problem: Franca-Hughes stabilization on cavity flow; Velocity error comparison between the \textit{offline-online stabilization} with/without supremizer and \textit{offline-only stabilization} with supremizer using $\mathbb{P}_2/{\mathbb{P}_2}$; stabilization coefficient $\delta=0.05, 0.5$; $N_u=N_p=N_s=20.$}
    \label{fig: U_alpha0_05_steady}
\end{figure}

\begin{figure}[H]
  \centering
  \begin{subfigure}{0.48\textwidth}
    \centering\includegraphics[width=\textwidth]{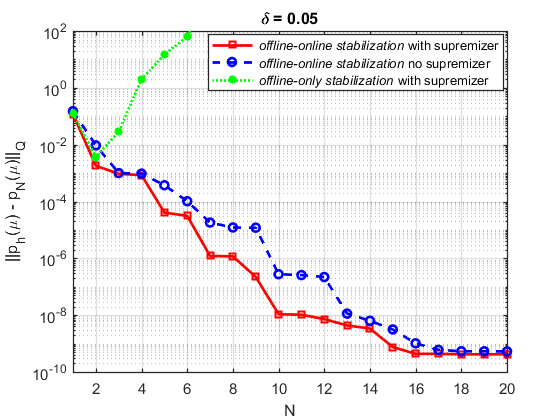}   
  \end{subfigure}%
  \quad%
  \begin{subfigure}{0.48\textwidth}
    \centering\includegraphics[width=\textwidth]{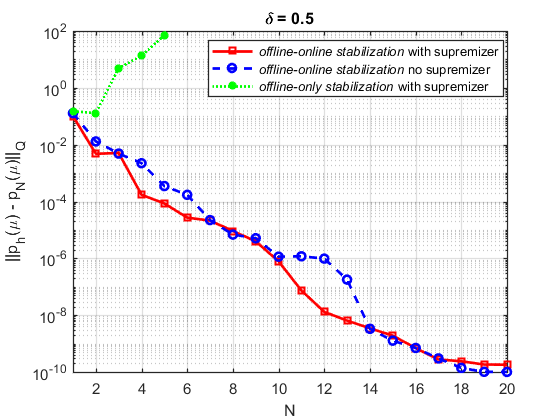}
  \end{subfigure}
     \caption{Stokes problem: Franca-Hughes stabilization on cavity flow; Pressure error comparison between the \textit{offline-online stabilization} with/without supremizer and \textit{offline-only stabilization} with supremizer using $\mathbb{P}_2/{\mathbb{P}_2}$; stabilization coefficient $\delta=0.05, 0.5$; $N_u=N_p=N_s=20.$}
    \label{fig: P_alpha0_05_steady}
\end{figure}

 
Figures \ref{fig: U_alpha0_05_steady} and \ref{fig: P_alpha0_05_steady} are plotted to see the error comparison for velocity and pressure, respectivley using Franca-Hughes stabilization \eqref{eq:our_choice} for $\mathbb{P}_2/\mathbb{P}_2$ FE pair and by varying the stabilization coefficient $\delta = 0.5, 0.05$. If we compare these results with previous results obtained for the $\mathbb{P}_1/\mathbb{P}_1$ case, we see that Franca-Hughes stabilization is able to perform comparably between the cases of \textit{offline-online stabilization} with and without supremizer, even for the pressure approximation. Therefore, we conclude that when we use Franca-Hughes stabilization with $\mathbb{P}_2/\mathbb{P}_2$ FE pair, there is no need to enrich the RB velocity space with supremizer solutions. In this way we can reduce also the online computational cost by decreasing the dimension of reduced velocity space (see Table \ref{tab:comput_steady_stokes_GEO}).
Furthermore, the results for different values of $\delta$ show the robustness of the methodology when varying the stabilization coefficient.


\begin{center}
 \begin{tabular}{|| l | l||}
 \hline\hline
Physical parameter & $\mu_1$ (fluid viscosity)\\ 
 \hline
Geometrical parameter & $\mu_2$ (horizontal length of domain)\\ 
 \hline
Range of $\mu_1$  & [0.25,0.75]  \\ 
 \hline
Range of $\mu_2$  & [1,3]  \\ 
 \hline
$\mu_1$ \textit{online}  & 0.6  \\ 
 \hline
$\mu_2$ \textit{online}  & 2  \\ 
 \hline
\multirow{2}{*}{FE degrees of freedom} & $6222$ ($\mathbb{P}_1/{\mathbb{P}_1}$)\\
& $10935$ ($\mathbb{P}_2/{\mathbb{P}_2}$)\\
 \hline
RB dimension & $N_u=N_s=N_p=20$\\
 \hline
Computation time ($\mathbb{P}_2/{\mathbb{P}_1}$) & $260s$ (\textit{offline}), $12s$ (\textit{online}) with supremizer\\
 \hline
\multirow{3}{*}{\textit{Offline} time ($\mathbb{P}_1/{\mathbb{P}_1}$)} & $180s$ (\textit{offline-online stabilization} with supremizer)\\
& $130s$ (\textit{offline-online stabilization} without supremizer)\\
& $105s$ (\textit{offline-only stabilization} with supremizer)\\
 \hline
\multirow{3}{*}{\textit{Offline} time ($\mathbb{P}_2/{\mathbb{P}_2}$)} & $348s$ (\textit{offline-online stabilization} with supremizer)\\
& $309s$ (\textit{offline-online stabilization} without supremizer)\\
& $280s$ (\textit{offline-only stabilization} with supremizer)\\
 \hline
\multirow{3}{*}{\textit{Online} time ($\mathbb{P}_1/{\mathbb{P}_1}$)} & $10s$ (\textit{offline-online stabilization} with supremizer)\\
& $8s$ (\textit{offline-online stabilization} without supremizer)\\
& $7s$ (\textit{offline-only stabilization} with supremizer)\\
 \hline
\multirow{3}{*}{\textit{Online} time ($\mathbb{P}_2/{\mathbb{P}_2}$)} & $15s$ (\textit{offline-online stabilization} with supremizer)\\
& $13s$ (\textit{offline-online stabilization} without supremizer)\\
& $10s$ (\textit{offline-only stabilization} with supremizer)\\
\hline\hline
\end{tabular}\vspace{0.3cm}
\captionof{table}{Stokes problem: Computational details for physical and geometrical parameters. }
\label{tab:comput_steady_stokes_GEO} 
\end{center}

From these results it is clear that there is a trade-off between increasing by supremizers reduced basis velocity spaces and improving accuracy of results.

\begin{figure}[H]
  \centering
  \begin{subfigure}{0.48\textwidth}
    \centering\includegraphics[width=\textwidth]{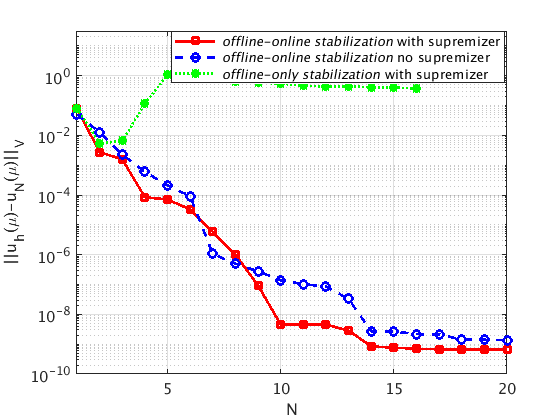}   
  \end{subfigure}%
  \quad%
  \begin{subfigure}{0.48\textwidth}
    \centering\includegraphics[width=\textwidth]{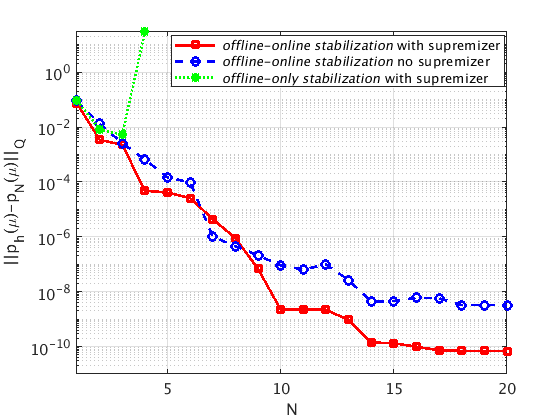}
  \end{subfigure}
  \caption{Stokes problem: Brezzi-Pitkaranta stabilization on cavity flow; Velocity (left) and pressure (right) error comparison between the \textit{offline-online stabilization} with/without supremizer and \textit{offline-only stabilization} with supremizer using $\mathbb{P}_1/{\mathbb{P}_1}$; stabilization coefficient $\delta=0.05$; $N_u=N_p=N_s=20.$}
\label{fig: U_alpha0_05_BP}
\end{figure}

\subsection{Numerical results for $\mathbb{P}_1/{\mathbb{P}_0}$}
\label{fig:stokes Sol_P1P0}
In this section we discuss the solution of steady parametrized Stokes problem using the FE pair given by $\mathbb{P}_1$ approximation for velocity and discontinuous $\mathbb{P}_0$ discretization for pressure. We present this case to show that the proposed methodology is applicable to any type of stabilization, and not necessarily residual based. Indeed, for $\mathbb{P}_1/{\mathbb{P}_0}$, both \eqref{eq:choice_general}-\eqref{eq:our_choice} vanish; one can alternatively resort to \cite{QV}:
\begin{equation}
s^{p,q}_{h}(q_{h};\boldsymbol\mu):=\delta\sum_{\sigma\in\Gamma_h}h_{\sigma}\int_{\sigma}\left[p_{h}\right]_{\sigma}\left[q_{h}\right]_{\sigma} 
\label{eq:P1P0 STAB}
\end{equation}
where $\Gamma_h$ is the set of all edges $\sigma$ of the triangulation except for those belonging to the boundary $\partial\Omega$, $h_{\sigma}$ is the length of $\sigma$ and $\left[q_{h}\right]_{\sigma}$ denotes its jump across $\sigma$.\par
\begin{figure}[H]
  \centering
  \begin{subfigure}{0.48\textwidth}
    \centering\includegraphics[width=\textwidth]{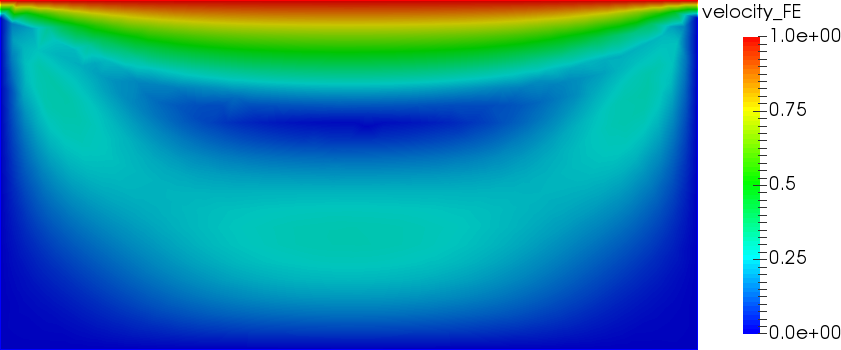}
    \caption{FE Velocity}
  \end{subfigure}%
  \quad%
  \begin{subfigure}{0.48\textwidth}
    \centering\includegraphics[width=\textwidth]{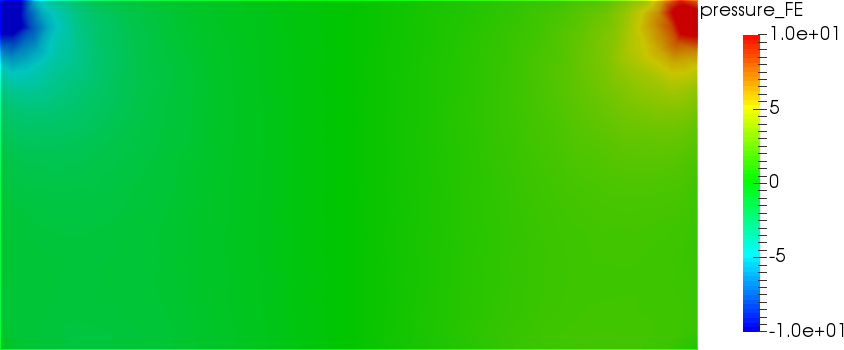}
    \caption{FE Pressure}
  \end{subfigure}
    \quad%
      \begin{subfigure}{0.48\textwidth}
    \centering\includegraphics[width=\textwidth]{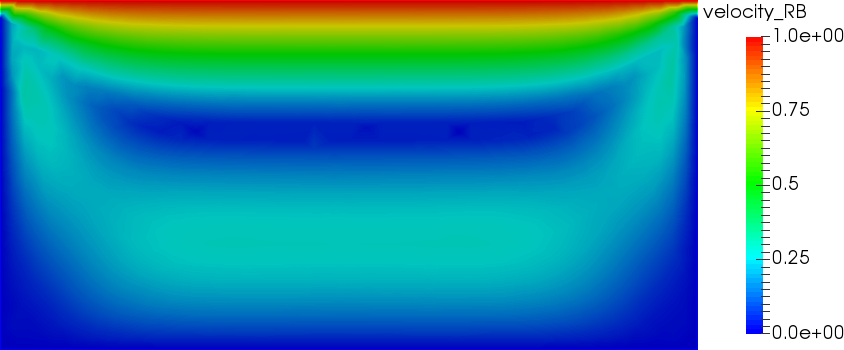}
    \caption{RB Velocity: \textit{offline-online stabilization} without supremizers}
  \end{subfigure}%
  \quad%
  \begin{subfigure}{0.48\textwidth}
    \centering\includegraphics[width=\textwidth]{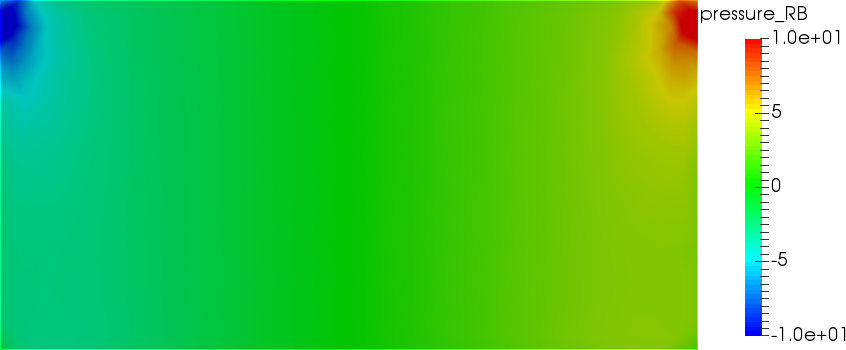}
    \caption{RB Pressure: \textit{offline-online stabilization} without supremizers}
  \end{subfigure}
  \quad%
   \begin{subfigure}{0.48\textwidth}
    \centering\includegraphics[width=\textwidth]{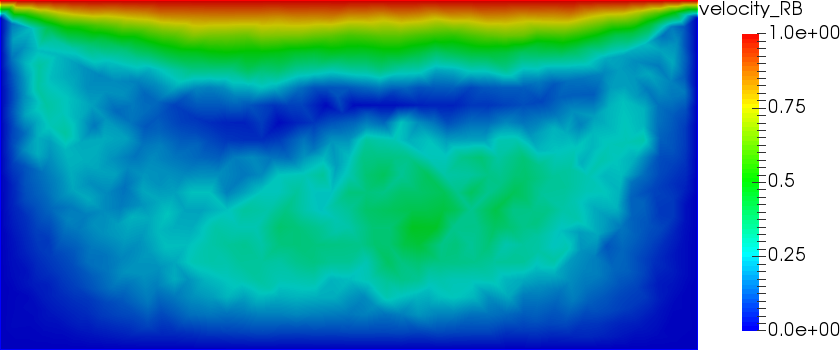}
    \caption{RB Velocity: \textit{offline-only stabilization} with supremizers}
  \end{subfigure}%
  \quad%
  \begin{subfigure}{0.48\textwidth}
    \centering\includegraphics[width=\textwidth]{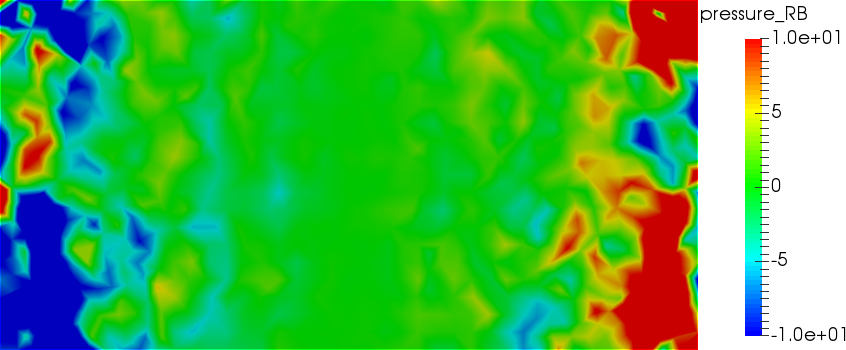}
    \caption{RB Pressure: \textit{offline-only stabilization} with supremizers }
  \end{subfigure}
    \quad%
  \caption{Stokes problem: FE and RB solutions for velocity and pressure at $(\mu_1,\mu_2)=(0.6,2)$; $N_u=N_p=20$ and using $\mathbb{P}_1/\mathbb{P}_0$.}
  \label{fig:P1P0}
\end{figure}

\begin{figure}[H]
  \centering
  \begin{subfigure}{0.48\textwidth}
    \centering\includegraphics[width=\textwidth]{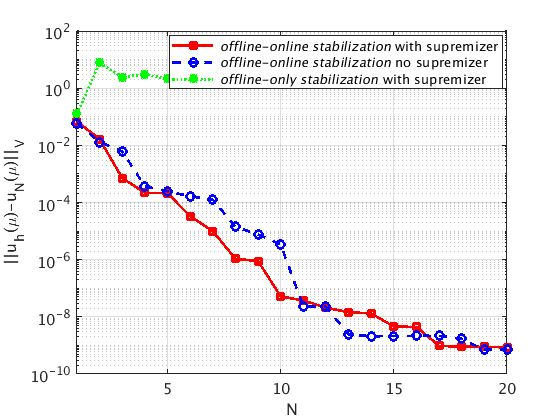}   
  \end{subfigure}%
  \quad%
  \begin{subfigure}{0.48\textwidth}
    \centering\includegraphics[width=\textwidth]{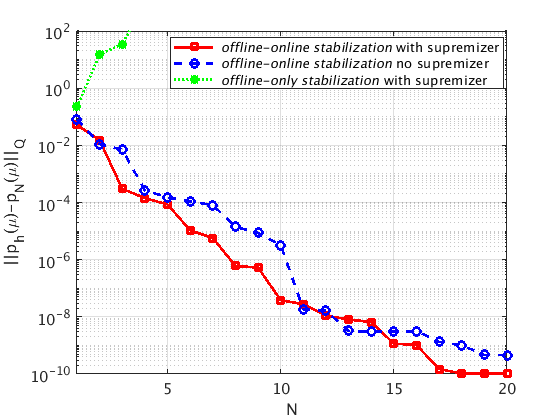}
  \end{subfigure}
  \caption{Stokes problem: Stabilization with $\mathbb{P}_1/{\mathbb{P}_0}$ on cavity flow: Velocity (left) and pressure (right) error between FE solution and RB solution for different possible options.}
\label{fig: U_P1P0}
\end{figure}

 
In Fig. \ref{fig:P1P0} we show some snapshots for velocity and pressure fields using \textit{stabilized} FE method and \textit{stabilized} RB method. From these solution plots we conclude that we are able to recover a good qualitative approximation of FE solution at reduced order level using the \textit{offline-online stabilization} without supremizers, whereas the \textit{offline-only stabilization} is not enough to recover FE approximation.
In Figs. \ref{fig: U_P1P0} we plot the comparison between \textit{offline-online stabilization} with/without supremizer and \textit{offline-only stabilization} with supremizer for velocity and pressure, respectively for $\mathbb{P}_1/{\mathbb{P}_0}$. These comparison shows that the \textit{offline-online stabilization} results in the most accurate method, and that the addition of supremizer to velocity space is not necessary for pressure recovery. Indeed, a good approximation of pressure is obtained even without the supremizer, and the enrichment accounts for a further improvement of only one order of magnitude. 

\label{subsec:StokesP1P0}

\subsection{Numerical results for stable $\mathbb{P}_2/{\mathbb{P}_1}$}
\begin{figure}[H]
  \centering
  \begin{subfigure}{0.48\textwidth}
    \centering\includegraphics[width=\textwidth]{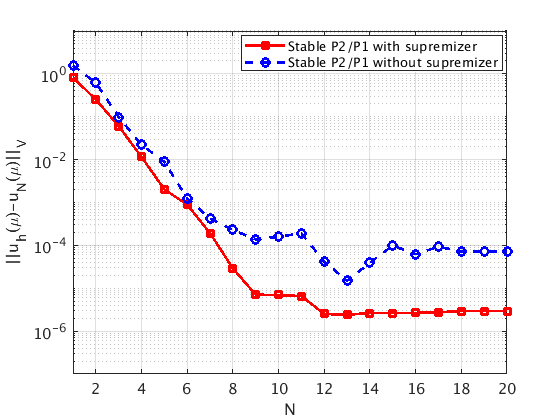}   
  \end{subfigure}%
  \quad%
  \begin{subfigure}{0.48\textwidth}
    \centering\includegraphics[width=\textwidth]{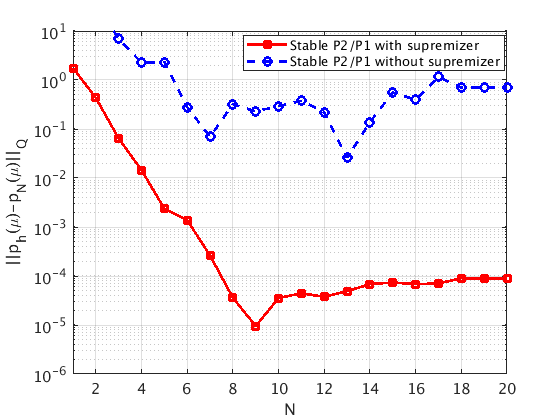}
  \end{subfigure}
 \caption{Stokes problem: Error between FE and RB solutions for velocity (left) and pressure (right) with/without supremizer using $\mathbb{P}_2/{\mathbb{P}_1}$.} 
\label{fig: UP_P2P1}\end{figure}

Just for completeness, here we plot a comparison between the FE solution and RB solution for velocity and pressure in Fig. \ref{fig: UP_P2P1} using stable FE pair $\mathbb{P}_2/\mathbb{P}_1$ \cite{Veroy2007}. We compare results with and without supremizer and conclude that supremizer is necessary to enrich the RB velocity space. In this case we are not using any \textit{stabilization} method, therefore we need supremizer to fulfill reduced inf-sup condition.
\label{subsec:StokesP2P1}
\pagebreak
\section{Parametrized Steady Navier-Stokes problem}
\label{sec: NS_continuous_formulation}
Let us now consider the steady incompressible Navier-Stokes equations in a parametrized domain $\Omega_o(\boldsymbol\mu)\subset\mathbb{R}^2$. The continuous parametrized formulation read as follows: find $\boldsymbol{u}_o(\boldsymbol\mu): \Omega_o(\boldsymbol\mu) \to \mathbb{R}^2$ and $p_o(\boldsymbol\mu): \Omega_o(\boldsymbol\mu) \to \mathbb{R}^2$ such that
\begin{equation}
\begin{cases}
-\nu(\boldsymbol\mu)\Delta{\boldsymbol{u}_o}(\boldsymbol\mu)+\left(\boldsymbol{u}_o(\boldsymbol\mu)\cdot\nabla\right)\boldsymbol{u}_o(\boldsymbol\mu)+\nabla{p}_o(\boldsymbol\mu)= \boldsymbol{0}& \text{ in } \Omega_o(\boldsymbol\mu),\  \\ 
\rm div\ {\boldsymbol{u}_o}(\boldsymbol\mu)=0 & \text{ in } \Omega_o(\boldsymbol\mu),\\
\boldsymbol{u}_o(\boldsymbol\mu)=\boldsymbol{g}_D(\boldsymbol\mu) & \text{ on } \Gamma_{D,o}(\boldsymbol\mu),\\
\boldsymbol{u}_o(\boldsymbol\mu)=\boldsymbol{0}, & \text{ on } \Gamma_{W,o}(\boldsymbol\mu),
\end{cases}
\label{eq:NS_continuous}
\end{equation}
where the parameter, unknowns, coefficients and boundaries have the same meaning as in the Stokes case \eqref{eq:S_continuous}. For the sake of the numerical test cases, we introduced the Reynolds number defined as $Re(\boldsymbol\mu)=L(\boldsymbol\mu)|\bar{\boldsymbol{u}}(\boldsymbol\mu)|/\nu(\boldsymbol\mu),$ being $L(\boldsymbol\mu)$ a characteristic length of the parametrized domain, $\bar{\boldsymbol{u}}(\boldsymbol\mu)$ a typical flow velocity and $\nu(\boldsymbol\mu)$ the kinematic viscosity. \par 
The weak formulation of problem \eqref{eq:NS_continuous} can be obtained, as in section 2, with the aid of a reference domain $\Omega$, on which function spaces $\boldsymbol{V}$ and $Q$ are defined, as in section 2. Proceeding as in the Stokes case, we have the following parametrized weak formulation of problem \eqref{eq:NS_continuous}:
\begin{equation}
\begin{cases}
\text{Find } \boldsymbol{u}(\boldsymbol\mu)\in \boldsymbol{V}, \, p(\boldsymbol\mu)\in Q:\\
a(\boldsymbol{u}(\boldsymbol\mu),\boldsymbol{v};\boldsymbol\mu)+b(\boldsymbol{v},p(\boldsymbol\mu);\boldsymbol\mu)+c(\boldsymbol{u}(\boldsymbol\mu),\boldsymbol{u}(\boldsymbol\mu),\boldsymbol{v};\boldsymbol\mu)+d(\boldsymbol{u}(\boldsymbol\mu),\boldsymbol{v};\boldsymbol\mu)=F(\boldsymbol{v};\boldsymbol\mu) & \forall \, v\in \boldsymbol{V},\\ b(\boldsymbol{u}(\boldsymbol\mu),q;\boldsymbol\mu)=G(q;\boldsymbol\mu) & \forall \, q\in Q,
\end{cases}
\label{eq:NS_weak}
\end{equation}
where $a(\boldsymbol{u},\boldsymbol{v};\boldsymbol\mu)$, $b(\boldsymbol{v},q;\boldsymbol\mu)$ are the bilinear forms related to diffusion and pressure-divergence operators, respectively defined in \eqref{eq: S_bilinear_forms}, whereas the trilinear form related to the convective term is defined as:
\begin{equation}
c(\boldsymbol{u},\boldsymbol{v},\boldsymbol{w};\boldsymbol\mu)=\int_{\Omega}u_i \chi_{ji}(x;\boldsymbol\mu)\dfrac{\partial{v_m}}{\partial x_j}w_m d\boldsymbol{x}.
\label{eq: NS_trilinear}
\end{equation} 
Other terms appearing due to the lifting of Dirichlet boundary conditions are defined as
\begin{equation}
\begin{split}
d(\boldsymbol{u},\boldsymbol{v};\boldsymbol\mu)&=c(\boldsymbol{l(\mu)},\boldsymbol{u},\boldsymbol{v};\boldsymbol\mu)+c(\boldsymbol{u},\boldsymbol{l(\mu)},\boldsymbol{v};\boldsymbol\mu),\\
F(\boldsymbol{v};\boldsymbol\mu)&=-a(\boldsymbol{l(\mu)},\boldsymbol{v},\boldsymbol\mu)-c(\boldsymbol{l(\mu)},\boldsymbol{l(\mu)},\boldsymbol{v};\boldsymbol\mu),\\
G(q;\boldsymbol\mu)&=-b(\boldsymbol{l(\mu)},q;\boldsymbol\mu),
\end{split}
\label{eq: lift_functions}
\end{equation}
being $\boldsymbol{l}(\boldsymbol\mu)$ a parametrized lifting function as in section 2.

\subsection{Finite Element formulation}
\label{sec: NS_FE_formulation}
In order to write the Galerkin-FE formulation for \eqref{eq:NS_weak}, we first need to introduce two finite-dimensional subspaces $\boldsymbol{V}_h\subset\boldsymbol{V}$, $Q_h\subset Q$ of dimension $\mathcal{N}_u$ and $\mathcal{N}_p$, respectively, as in section 2.1. The Galerkin-FE approximation of the parametrized problem \eqref{eq:NS_weak} reads as follows:
\begin{equation}
\begin{cases}
\text{Find } \boldsymbol{u}_{h}(\boldsymbol{\mu})\in \boldsymbol{V}_{h}, p_{h}(\boldsymbol{\mu})\in Q_{h}:\\
a(\boldsymbol{u}_h(\boldsymbol\mu),\boldsymbol{v}_h;\boldsymbol\mu)+b(\boldsymbol{v}_h,p_h(\boldsymbol\mu);\boldsymbol\mu)+c(\boldsymbol{u}_h(\boldsymbol\mu),\boldsymbol{u}_h(\boldsymbol\mu),\boldsymbol{v}_h;\boldsymbol\mu)\\+d(\boldsymbol{u}_h(\boldsymbol\mu),\boldsymbol{v}_h;\boldsymbol\mu)=F(\boldsymbol{v}_h;\boldsymbol\mu) & \forall \, \boldsymbol{v}_h\in \boldsymbol{V}_h,\\ b(\boldsymbol{u}_h(\boldsymbol\mu),q_h;\boldsymbol\mu)=G(q_h;\boldsymbol\mu) & \forall \, q_h\in Q_h,
\end{cases}
\label{eq:NS_weak_FE}
\end{equation}
where
\begin{equation}
\begin{split}
d(\boldsymbol{u}_h(\boldsymbol\mu),\boldsymbol{v}_h;\boldsymbol\mu)&=c(\boldsymbol{l}_h(\boldsymbol\mu),\boldsymbol{u}_h,\boldsymbol{v}_h;\boldsymbol\mu)+c(\boldsymbol{u}_h,\boldsymbol{l}_h(\boldsymbol\mu),\boldsymbol{v}_h;\boldsymbol\mu),\\
F(\boldsymbol{v}_h;\boldsymbol\mu)&=-a(\boldsymbol{l}_h(\boldsymbol\mu),\boldsymbol{v}_h,\boldsymbol\mu)-c(\boldsymbol{l}_h(\boldsymbol\mu),\boldsymbol{l}_h(\boldsymbol\mu),\boldsymbol{v}_h;\boldsymbol\mu),\\
G(q_h;\boldsymbol\mu)&=-b(\boldsymbol{l}_h(\boldsymbol\mu),q_h;\boldsymbol\mu).
\end{split}
\label{eq: lift_functions_FE}
\end{equation}
Let $\{\boldsymbol\phi_i^h\}_{i=1}^{\mathcal{N}_u}$ and $\{\psi_j^h\}_{j=1}^{\mathcal{N}_p}$ be basis functions of $\boldsymbol{V}_h$ and $Q_h$ respectively. We introduce the matrices $A(\boldsymbol\mu)\in\mathbb{R}^{\mathcal{N}_u\times\mathcal{N}_u}$,
$C(\boldsymbol{u}(\boldsymbol\mu);\boldsymbol\mu)\in\mathbb{R}^{\mathcal{N}_u\times\mathcal{N}_u}$, and $B(\boldsymbol\mu)\in\mathbb{R}^{\mathcal{N}_p\times\mathcal{N}_u}$ whose entries are
\begin{equation}
\begin{split}
\left(A(\boldsymbol\mu)\right)_{ij}&=a(\boldsymbol\phi_j^h,\boldsymbol\phi_i^h;\boldsymbol\mu)+d(\boldsymbol\phi_j^h,\boldsymbol\phi_i^h;\boldsymbol\mu), \qquad
\left(B(\boldsymbol\mu)\right)_{ki}=b(\boldsymbol\phi_i^h,\boldsymbol\psi_k^h;\boldsymbol\mu),\\
\left(C(\boldsymbol{u}(\boldsymbol\mu);\boldsymbol\mu)\right)_{ij}&=\sum_{m=1}^{\mathcal{N}_u}\boldsymbol{u}_h^{m}(\boldsymbol\mu)c(\boldsymbol\phi_m^h,\boldsymbol\phi_j^h,\boldsymbol\phi_i^h;\boldsymbol\mu), \quad \text{ for } 1\leq i,j\leq\mathcal{N}_u, 1\leq k\leq\mathcal{N}_p,
\end{split}
\label{eq:LHS_matrices}
\end{equation}
and therefore the nonlinear algebraic system is:
\begin{equation}
\left[\begin{array}{cc} A(\boldsymbol\mu)+C(\boldsymbol{u}_h(\boldsymbol\mu);\boldsymbol\mu) & B^T(\boldsymbol\mu)\\ B(\boldsymbol\mu) & \boldsymbol{0} \end{array}\right]
\left[\begin{array}{c} \boldsymbol{U}(\boldsymbol\mu) \\ \boldsymbol{P}(\boldsymbol\mu) \end{array}\right] 
= \left[\begin{array}{c} \boldsymbol{\bar{f}}(\boldsymbol\mu) \\ \boldsymbol{\bar{g}}(\boldsymbol\mu) \end{array}\right], 
\label{eq: Algebraic_Navier_Stokes}
\end{equation}
for the vectors of coefficients  $\boldsymbol{U}=(u_h^{(1)},...,u_h^{({\mathcal{N}_u})})^T, \boldsymbol{P}=(p_h^{(1)},...,p_h^{({\mathcal{N}_p})})^T$ such that $\boldsymbol{u}_h = \sum_{i=1}^{\mathcal{N}_u} u_h^{(i)} \boldsymbol\phi_i^h$ and $p_h = \sum_{j=1}^{\mathcal{N}_p} p_h^{(j)} \psi_j^h$, where for $1\leq i\leq\mathcal{N}_u$ and $1\leq k\leq\mathcal{N}_p$:
\begin{equation}
(\boldsymbol{\bar{f}}(\boldsymbol\mu))_i=-a(\boldsymbol{l}_h,\boldsymbol\phi_i^h;\boldsymbol\mu)-c(\boldsymbol{l}_h,\boldsymbol{l}_h,\boldsymbol\phi_i^h;\boldsymbol\mu), \quad (\boldsymbol{\bar{g}}(\boldsymbol\mu))_k=-b(\boldsymbol{l}_h,\boldsymbol\psi_k^h;\boldsymbol\mu), 
\label{eq:RHS_matrices}
\end{equation}
with $\boldsymbol{l}_h=\boldsymbol{l}_h(\boldsymbol\mu),$ the FE interpolant of lifting function. 
To solve this nonlinear system we use the Newton method \cite{Quarteroni2000}. \par
For an efficient RB method, we need to ensure the assumption of affine parametric dependence on operators \eqref{eq:LHS_matrices} and \eqref{eq:RHS_matrices}, i.e, these operators can be written as:
\begin{equation}
\begin{split}
A(\boldsymbol\mu)&=\sum_{q=1}^{Q_a}\Theta_{q}^{a}(\boldsymbol\mu)A^q, \qquad C(.;\boldsymbol\mu)=\sum_{q=1}^{Q_c}\Theta_{q}^{c}(\boldsymbol\mu)C^q(.), \qquad B(\boldsymbol\mu)=\sum_{q=1}^{Q_b}\Theta_{q}^{b}(\boldsymbol\mu)B^{q}, \\
\boldsymbol{\bar{f}}(\boldsymbol\mu)&=\sum_{q=1}^{Q_f}\Theta_{q}^{f}(\boldsymbol\mu)\boldsymbol{\bar{f}}^q,\qquad  \boldsymbol{\bar{g}}(\boldsymbol\mu)=\sum_{q=1}^{Q_g}\Theta_{q}^{g}(\boldsymbol\mu)\boldsymbol{\bar{g}}^q.
\end{split}
\label{eq: affine_NS}
\end{equation}
\subsection{Stabilized Finite Element formulation}
\label{sec: NS_FE_stabilization}
In this section we introduce the stabilization terms into the FE formulation of \eqref{eq:NS_weak_FE}. The stabilized FE formulation of \eqref{eq:NS_weak_FE} read as:
\begin{equation}
\begin{cases}
\text{Find } \boldsymbol{u}_{h}(\boldsymbol{\mu})\in \boldsymbol{V}_{h}, p_{h}(\boldsymbol{\mu})\in Q_{h}:\\ 
a(\boldsymbol{u}_h(\boldsymbol\mu),\boldsymbol{v}_h;\boldsymbol\mu)+b(\boldsymbol{v}_h,p_h(\boldsymbol\mu);\boldsymbol\mu)+c(\boldsymbol{u}_h(\boldsymbol\mu),\boldsymbol{u}_h(\boldsymbol\mu),\boldsymbol{v}_h;\boldsymbol\mu)\\+d(\boldsymbol{u}_h(\boldsymbol\mu),\boldsymbol{v}_h;\boldsymbol\mu)-s^{u,v}_{h}(\boldsymbol{u}_{h},\boldsymbol{v}_{h};\boldsymbol\mu)-s^{p,v}_{h}(p_h,\boldsymbol{v}_{h};\boldsymbol\mu)=F(\boldsymbol{v}_h;\boldsymbol\mu) & \forall \, \boldsymbol{v}_h\in \boldsymbol{V}_h,\\ b(\boldsymbol{u}_h(\boldsymbol\mu),q_h;\boldsymbol\mu)-s^{u,q}_{h}(\boldsymbol{u}_{h},q_{h};\boldsymbol\mu)-s^{p,q}_{h}(p_h,q_{h};\boldsymbol\mu)=G(q_h;\boldsymbol\mu) & \forall \, q_h\in Q_h,
\end{cases}
\label{eq:NS_weak_stab_FE}
\end{equation}
where $s^{u,v}_{h}(.,.;\boldsymbol\mu)$, $s^{p,v}_{h}(.,.;\boldsymbol\mu)$, $s^{u,q}_{h}(.,.;\boldsymbol\mu)$ and $s^{p,q}_{h}(.,.;\boldsymbol\mu)$ are the stabilization terms \cite{QV} defined as:
\begin{equation}
\begin{split}
s^{u,v}_{h}(\boldsymbol{u}_{h},\boldsymbol{v}_{h};\boldsymbol\mu)&:=\delta\sum_{K}h_K^{2}\int_K(-\nu\Delta\boldsymbol{u}_{h}+\boldsymbol{u}_h\cdot\nabla{\boldsymbol{u}_h},-\rho\nu\Delta\boldsymbol{v}_{h}),\\
s^{p,v}_{h}(p_h,\boldsymbol{v}_{h};\boldsymbol\mu)&:=\delta\sum_{K}h_K^{2}\int_K(\nabla{p}_{h},-\rho\nu\Delta\boldsymbol{v}_{h}),\\
s^{u,q}_{h}(\boldsymbol{u}_{h},q_{h};\boldsymbol\mu)&:=\delta\sum_{K}h_K^{2}\int_K(-\nu\Delta\boldsymbol{u}_{h}+\boldsymbol{u}_h\cdot\nabla{\boldsymbol{u}_h},\nabla{q}_{h}),\\
s^{p,q}_{h}(p_h,q_{h};\boldsymbol\mu)&:=\delta\sum_{K}h_K^{2}\int_K(\nabla{p}_{h},\nabla{q}_{h}),
\end{split}
\label{eq: stab_terms}
\end{equation}
where $\delta$ is the stabilization coefficient needs to be chosen properly \cite{Lube2006, Braack2007, Lube2006(1)}. For $\rho=0, 1$, the stabilization \eqref{eq: stab_terms} is respectively known as Streamline Upwind Petrov Galerkin (SUPG) \cite{Brooks1982}, Galerkin least-squares (GLS) \cite{Hughes1989}. The case $\rho=-1$ was studied by Franca and Frey \cite{FRANCA1992}.
Several other works on these kind of stabilization techniques can be found in \cite{Lube2002, Tobiska1996, ZHU1993, Lube1990} and references therein.
\par 
In this paper we discuss only the SUPG stabilization. Therefore, after adding the SUPG stabilization terms, the stabilized version of nonlinear system reads as \eqref{eq: Algebraic_Navier_Stokes}:
\begin{equation}
\left[\begin{array}{cc} A(\boldsymbol\mu)+\tilde{C}(\boldsymbol{u}_h(\boldsymbol\mu);\boldsymbol\mu) & \tilde{B^T}(\boldsymbol\mu)\\ \tilde{B}(\boldsymbol\mu) & -S_N(\boldsymbol\mu) \end{array}\right]
\left[\begin{array}{c} \boldsymbol{U}(\boldsymbol\mu) \\ \boldsymbol{P}(\boldsymbol\mu) \end{array}\right] 
= \left[\begin{array}{c} \boldsymbol{\bar{f}}(\boldsymbol\mu) \\ \boldsymbol{\bar{g}}(\boldsymbol\mu) \end{array}\right], 
\label{eq: Algebraic_Navier_Stokes_stab}
\end{equation}
where $\tilde{B},$ includes the the effects of stabilization in mass equation on divergence term, $\tilde{B^T}$ includes the effects on pressure gradient and $\tilde{C}$ contains the nonlinear stabilization terms \cite{QV}. 
\subsection{Reduced Basis formulation}
\label{subsec:RB_NS}
In this section we present the RB formulation of parametrized Navier-Stokes problem \eqref{eq:NS_continuous} \cite{quarteroni2007,Lovgren2006} in a similar way as we have done for the Stokes case. The RB spaces for velocity and pressure are defined in section \ref{subsec:RBStokes}.
\par 
Once, we have built the RB for velocity and pressure fields during the \textit{offline} stage, the RB formulation corresponding to FE formulation \eqref{eq:NS_weak_FE} reads as:
\begin{equation}
\begin{cases}
\text{Find } \boldsymbol{u}_N(\boldsymbol{\mu}) \in \tilde{\boldsymbol{V}}_N, p_N(\boldsymbol{\mu})) \in Q_N:\\
a(\boldsymbol{u}_N(\boldsymbol\mu),\boldsymbol{v}_N;\boldsymbol\mu)+b(\boldsymbol{v}_N,p_N(\boldsymbol\mu);\boldsymbol\mu)+c(\boldsymbol{u}_N(\boldsymbol\mu),\boldsymbol{u}_N(\boldsymbol\mu),\boldsymbol{v}_N;\boldsymbol\mu)\\+d(\boldsymbol{u}_N(\boldsymbol\mu),\boldsymbol{v}_N;\boldsymbol\mu)=F(\boldsymbol{v}_N;\boldsymbol\mu) & \forall \, \boldsymbol{v}_N\in \tilde{\boldsymbol{V}}_N ,\\ b(\boldsymbol{u}_N(\boldsymbol\mu),q_N;\boldsymbol\mu)=G(q_N;\boldsymbol\mu) & \forall \, q_N\in Q_N.
\end{cases}
\label{eq:NS_weak_RB}
\end{equation}
Owing to the representation \eqref{eq:RBsolution}, the corresponding nonlinear reduced system
\begin{equation}
\left[\begin{array}{cc} A_N(\boldsymbol\mu)+C_N(\boldsymbol{u}(\boldsymbol\mu);\boldsymbol\mu) & B_N^T(\boldsymbol\mu)\\ B_N(\boldsymbol\mu) & \boldsymbol{0} \end{array}\right]
\left[\begin{array}{c} \boldsymbol{U}_N(\boldsymbol\mu) \\ \boldsymbol{P}_N(\boldsymbol\mu) \end{array}\right] 
= \left[\begin{array}{c} \boldsymbol{\bar{f}}_N(\boldsymbol\mu) \\ \boldsymbol{\bar{g}}_N(\boldsymbol\mu) \end{array}\right], 
\label{eq: Algebraic_Navier_Stokes_RB}
\end{equation}
where the reduced order matrices are defined as:
\begin{equation}
\begin{split}
A_N(\boldsymbol\mu)&=Z^T_{u,s}A(\boldsymbol\mu)Z_{u,s}, \quad B_N(\boldsymbol\mu)=Z^T_{p}B(\boldsymbol\mu)Z_{u,s}, \quad C_N(.;\boldsymbol\mu)=Z^T_{u,s}C(.;\boldsymbol\mu)Z_{u,s},\\ \boldsymbol{\bar{f}}_N(\boldsymbol\mu)&=Z^T_{u,s}\boldsymbol{\bar{f}}(\boldsymbol\mu),\quad
\boldsymbol{\bar{g}}_N(\boldsymbol\mu)=Z^T_{p}\boldsymbol{\bar{g}}(\boldsymbol\mu),
\end{split}
\end{equation}
with $Z_{u,s}$, the velocity basis matrix and $Z_{p}$ denotes the pressure basis matrix.

\subsection{Stabilized Reduced Basis formulation}
In this section, we present the stabilized RB method for the model problem \eqref{eq:NS_continuous}. The stabilized RB formulation corresponding to stabilized FE formulation \eqref{eq:NS_weak_stab_FE} reads as:
\begin{equation}
\begin{cases}
\text{Find } \boldsymbol{u}_N(\boldsymbol{\mu}) \in \boldsymbol{V}_N, p_N(\boldsymbol{\mu})\in Q_N:\\
a(\boldsymbol{u}_N(\boldsymbol\mu),\boldsymbol{v}_N;\boldsymbol\mu)+b(\boldsymbol{v}_N,p_N(\boldsymbol\mu);\boldsymbol\mu)+c(\boldsymbol{u}_N(\boldsymbol\mu),\boldsymbol{u}_N(\boldsymbol\mu),\boldsymbol{v}_N;\boldsymbol\mu)\\+d(\boldsymbol{u}_N(\boldsymbol\mu),\boldsymbol{v}_N;\boldsymbol\mu)-s^{u,v}_{N}(\boldsymbol{u}_{N},\boldsymbol{v}_{N};\boldsymbol\mu)-s^{p,v}_{N}(p_N,\boldsymbol{v}_{N};\boldsymbol\mu)=F(\boldsymbol{v}_N;\boldsymbol\mu) & \forall \, \boldsymbol{v}_N\in \boldsymbol{V}_N,\\ b(\boldsymbol{u}_N(\boldsymbol\mu),q_N;\boldsymbol\mu)-s^{u,q}_{N}(\boldsymbol{u}_{N},q_{N};\boldsymbol\mu)-s^{p,q}_{N}(p_N,q_{N};\boldsymbol\mu)=G(q_N;\boldsymbol\mu) & \forall \, q_N\in Q_N,
\end{cases}
\label{eq:NS_weak_stab_RB}
\end{equation}
where $s^{u,v}_{N}(.,.;\boldsymbol\mu)$, $s^{p,v}_{N}(.,.;\boldsymbol\mu)$, $s^{u,q}_{N}(.,.;\boldsymbol\mu)$ and $s^{p,q}_{N}(.,.;\boldsymbol\mu)$ are the reduced order stabilization terms defined as:
\begin{equation}
\begin{split}
s^{u,v}_{N}(\boldsymbol{u}_{N},\boldsymbol{v}_{N};\boldsymbol\mu)&:=\delta\sum_{K}h_K^{2}\int_K(-\nu\Delta\boldsymbol{u}_{N+\boldsymbol{u}_N\cdot\nabla{\boldsymbol{u}_N}},-\rho\nu\Delta\boldsymbol{v}_{N}),\\
s^{p,v}_{N}(p_h,\boldsymbol{v}_{N};\boldsymbol\mu)&:=\delta\sum_{K}h_K^{2}\int_K(\nabla{p}_{N},-\rho\nu\Delta\boldsymbol{v}_{N}),\\
s^{u,q}_{N}(\boldsymbol{u}_{N},q_{N};\boldsymbol\mu)&:=\delta\sum_{K}h_K^{2}\int_K(-\nu\Delta\boldsymbol{u}_{N}+\boldsymbol{u}_N\cdot\nabla{\boldsymbol{u}_N},\nabla{q}_{N}),\\
s^{p,q}_{N}(p_N,q_N;\boldsymbol\mu)&:=\delta\sum_{K}h_K^{2}\int_K(\nabla{p}_{N},\nabla{q}_{N}),
\end{split}
\label{eq: stab_terms_RB}
\end{equation}
The corresponding nonlinear system results in
\begin{equation}
\left[\begin{array}{cc} A_N(\boldsymbol\mu)+\tilde{C}_N(\boldsymbol{u}(\boldsymbol\mu);\boldsymbol\mu) & \tilde{B}^T_N(\boldsymbol\mu)\\ \tilde{B}_N(\boldsymbol\mu) & -S_N(\boldsymbol\mu) \end{array}\right]
\left[\begin{array}{c} \boldsymbol{U}_N(\boldsymbol\mu) \\ \boldsymbol{P}_N(\boldsymbol\mu) \end{array}\right] 
= \left[\begin{array}{c} \boldsymbol{\bar{f}}_N(\boldsymbol\mu) \\ \boldsymbol{\bar{g}}_N(\boldsymbol\mu) \end{array}\right], 
\label{eq: Algebraic_Navier_Stokes_stab_RB}
\end{equation}
where $\tilde{B}_N, \tilde{B}^T_N$ and $\tilde{C}_N$ are RB stabilization matrices \cite{QV} defined as:
\begin{equation}
\tilde{B}_N(\boldsymbol\mu)=Z^T_{p}\tilde{B}(\boldsymbol\mu)Z_{u,s}, \quad \tilde{B}^{T}_N(\boldsymbol\mu)=Z_{p}\tilde{B}(\boldsymbol\mu)Z_{u,s}, \quad S_N(\boldsymbol\mu)=Z^T_{p}S(\boldsymbol\mu)Z_{p}.
\label{eq:RBmatrices_NavierStokes_STAB}
\end{equation}

\label{subsec: NS_RB_stabilization}


\section{Numerical results and discussion for parametrized Navier-Stokes problems}
\label{sec:Results_NS}
In this section we present some numerical results for the RB approximation of steady parametrized Navier-Stokes problem developed in section \ref{sec: NS_continuous_formulation} and subsections therein. As in section 3, we compare three possible options $(i)$ \textit{offline-online stabilization} with supremizer, $(ii)$ \textit{offline-online stabilization} without supremizer, $(iii)$ \textit{offline-only stabilization} with supremizer. The numerical test case is the same parametrized cavity problem as in section 3. Table \ref{tab:comput_steady_NS_GEO} illustrates the computation details and the cost of stabilization options using both $\mathbb{P}_1/{\mathbb{P}_1}$ and $\mathbb{P}_2/{\mathbb{P}_2}$ FE pairs. 
\par
\begin{figure}[H]
  \centering
  \begin{subfigure}{0.44\textwidth}
    \centering\includegraphics[width=\textwidth]{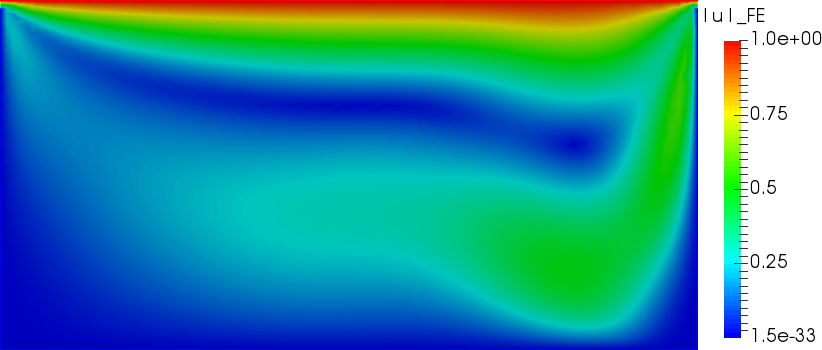}
    \caption{FE Velocity}
  \end{subfigure}%
  \quad%
  \begin{subfigure}{0.44\textwidth}
    \centering\includegraphics[width=\textwidth]{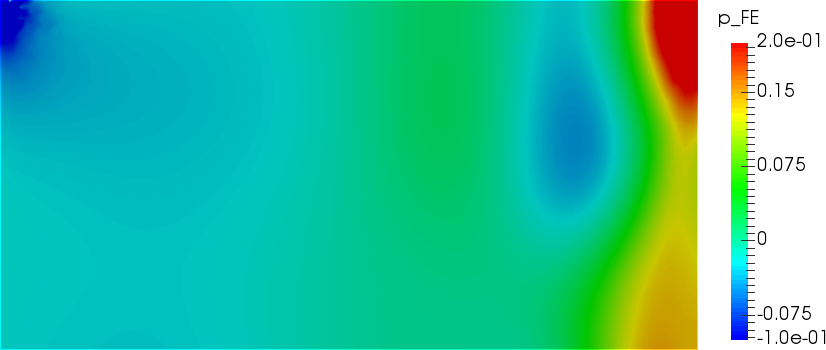}
    \caption{FE Pressure}
  \end{subfigure}
    \quad%
   \begin{subfigure}{0.44\textwidth}
    \centering\includegraphics[width=\textwidth]{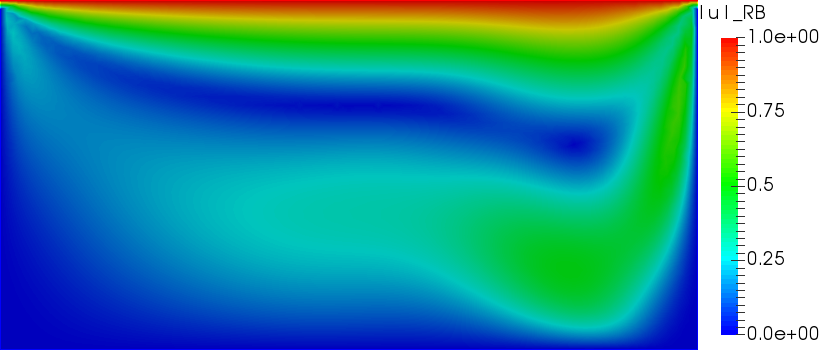}
    \caption{RB Velocity: \textit{offline-online}}
  \end{subfigure}%
  \quad%
  \begin{subfigure}{0.44\textwidth}
    \centering\includegraphics[width=\textwidth]{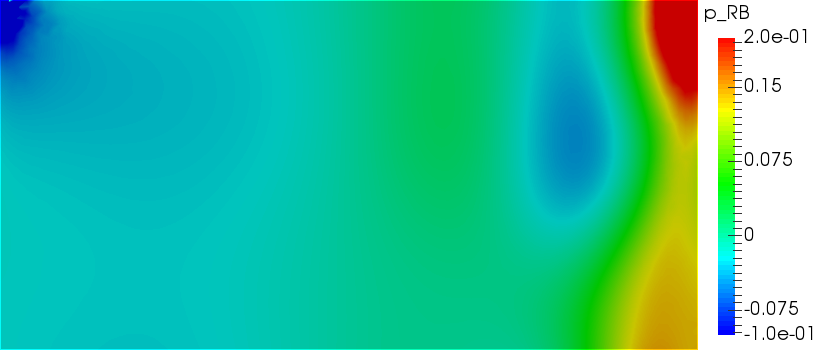}
    \caption{RB Pressure: \textit{offline-online}}
  \end{subfigure}
    \quad%
    \begin{subfigure}{0.44\textwidth}
    \centering\includegraphics[width=\textwidth]{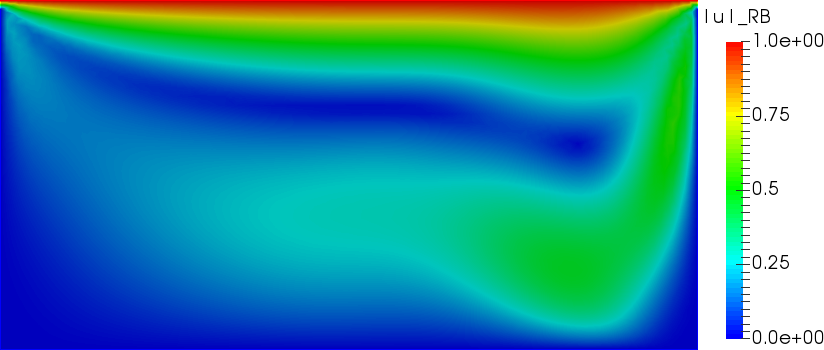}
    \caption{RB Velocity: \textit{offline-only}}
  \end{subfigure}%
  \quad%
  \begin{subfigure}{0.44\textwidth}
    \centering\includegraphics[width=\textwidth]{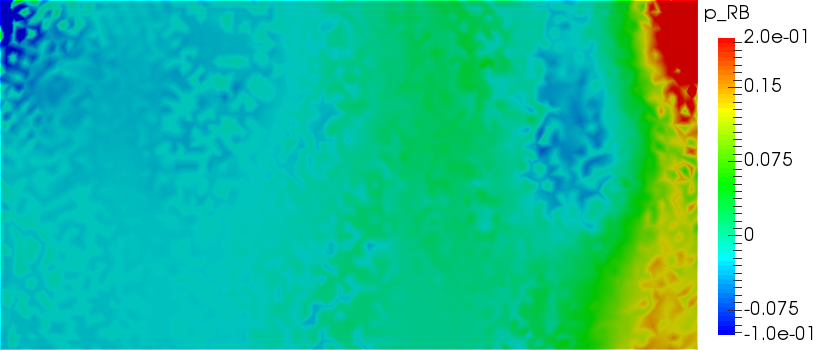}
    \caption{RB Pressure: \textit{offline-only}}
  \end{subfigure}
  \quad%
  \caption{Navier-Stokes problem with SUPG stabilization: From top to bottom; FE solution (top), RB solution using \textit{offline-online stabilization} (center), RB solution using \textit{offline-only stabilization} (bottom) for $(\mu_1,\mu_2)=(120,2).$}
\label{fig:SUPGEO_Sol}
\end{figure}

 
In Fig. \ref{fig:SUPGEO_Sol}, we show the FE solution (top) for velocity and pressure, RB solution for velocity and pressure using the \textit{offline-online stabilization} with/without supremizer (center) and the \textit{offline-only stabilization} with supremizer (bottom), for $(\mu_1,\mu_2)=(120,2)$ and $\delta=1.0$. These results are carried out using $\mathbb{P}_2/{\mathbb{P}_2}$ FE pair (results are similar for $\mathbb{P}_1/{\mathbb{P}_1}$). From these results, we see that the FE solution and RB solution obtained by \textit{offline-online stabilization} are same, but the RB pressure obtained by \textit{offline-only stabilization} is not accurate. 
\par
\begin{figure}[H]
  \centering
  \begin{subfigure}{0.48\textwidth}
    \centering\includegraphics[width=\textwidth]{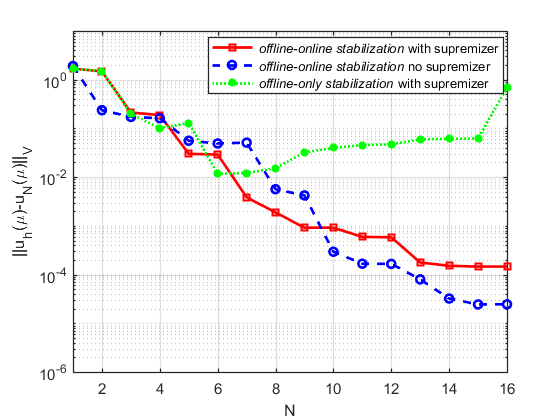}   
  \end{subfigure}%
  \quad%
  \begin{subfigure}{0.48\textwidth}
    \centering\includegraphics[width=\textwidth]{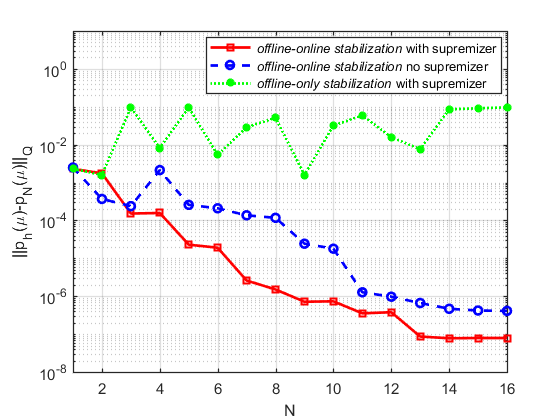}
  \end{subfigure}
  \caption{Navier-Stokes problem: SUPG stabilization with geometrical and physical parametrization on cavity flow; Error between FE and RB solution for velocity (left) and pressure (right) using $\mathbb{P}_1/{\mathbb{P}_1}$.}
    \label{fig:USUPG}
\end{figure}

\begin{figure}[H]
  \centering
  \begin{subfigure}{0.48\textwidth}
    \centering\includegraphics[width=\textwidth]{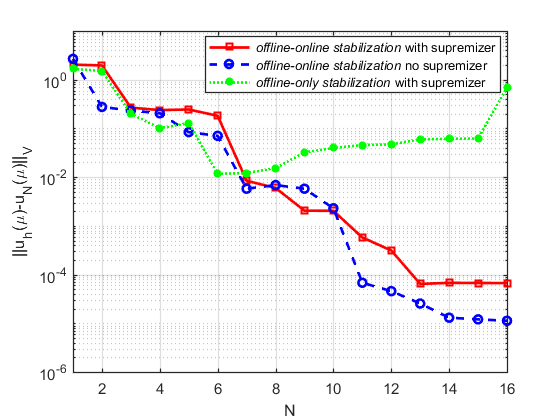}   
  \end{subfigure}%
  \quad%
  \begin{subfigure}{0.48\textwidth}
    \centering\includegraphics[width=\textwidth]{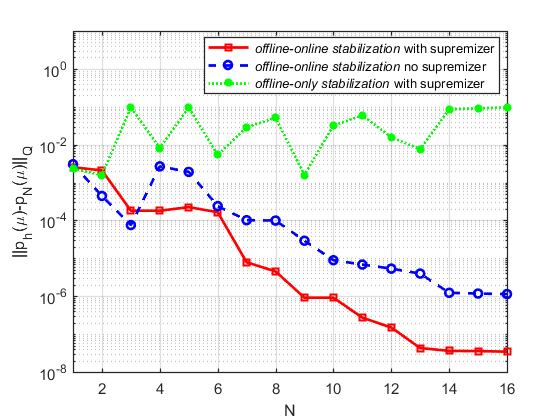}
  \end{subfigure}
  \caption{Navier-Stokes problem: SUPG stabilization with geometrical and physical parametrization on cavity flow; Error between FE and RB solution for velocity (left) and pressure (right) using $\mathbb{P}_2/{\mathbb{P}_2}$.}
    \label{fig:USUPG_P2}
\end{figure}

 
Figures \ref{fig:USUPG} and \ref{fig:USUPG_P2} represents the error between FE and RB velocity obtained by three possible stabilization options using $\mathbb{P}_1/{\mathbb{P}_1}$ and $\mathbb{P}_2/{\mathbb{P}_2}$ FE pair, respectively. From these results we conclude that the \textit{offline-online stabilization} without supremizer has better performance in case of velocity, when compared to other two options and is the most consistent stabilization option. However, in case of pressure, supremizer is improving the error upto an order of magnitude. 
Moreover, when comparing the computational costs of the two \textit{offline-online stabilization options, the one without supremizer has less computation time than the case with supremizer.} Therefore, we conclude that the \textit{offline-online stabilization} is sufficient to guarantee a stable RB solution and we can avoid the supremizer enrichment to increase the computational performance of \textit{online} stage. 
\par 

\begin{center}
 \begin{tabular}{|| l | l||}
 \hline\hline
Physical parameter & $\mu_1$ (Reynolds number)\\ 
 \hline
Geometrical parameter & $\mu_2$ (horizontal length of domain)\\ 
 \hline
Range of $\mu_1$  & [100,200]  \\ 
 \hline
Range of $\mu_2$  & [1.5,3]  \\ 
 \hline
$\mu_1$ \textit{online}  & 120  \\ 
 \hline
$\mu_2$ \textit{online}  & 2  \\ 
 \hline
\multirow{2}{*}{FE degrees of freedom} & $11160$ ($\mathbb{P}_1/{\mathbb{P}_1}$)\\
& $44091$ ($\mathbb{P}_2/{\mathbb{P}_2}$)\\
 \hline
RB dimension & $N_u=N_s=N_p=16$\\
 \hline
Computation time ($\mathbb{P}_2/{\mathbb{P}_1}$) & $3909s$ (\textit{offline}), $195s$ (\textit{online}) with supremizers\\
 \hline
\multirow{3}{*}{\textit{Offline} time ($\mathbb{P}_1/{\mathbb{P}_1}$)} & $2034s$ (\textit{offline-online stabilization} with supremizers)\\
& $1920s$ (\textit{offline-online stabilization} without supremizers)\\
& $649$ (\textit{offline-only stabilization} with supremizers)\\
 \hline
\multirow{3}{*}{\textit{Offline} time ($\mathbb{P}_2/{\mathbb{P}_2}$)} & $4885s$ (\textit{offline-online stabilization} with supremizers)\\
& $4387ss$ (\textit{offline-online stabilization} without supremizers)\\
& $1650s$ (\textit{offline-only stabilization} with supremizers)\\
 \hline
\multirow{3}{*}{\textit{Online} time ($\mathbb{P}_1/{\mathbb{P}_1}$)} & $110s$ (\textit{offline-online stabilization} with supremizers)\\
& $87s$ (\textit{offline-online stabilization} without supremizers)\\
& $35s$ (\textit{offline-only stabilization} with supremizers)\\
 \hline
\multirow{3}{*}{\textit{Online} time ($\mathbb{P}_2/{\mathbb{P}_2}$)} & $242s$ (\textit{offline-online stabilization} with supremizers)\\
& $180s$ (\textit{offline-online stabilization} without supremizers)\\
& $90s$ (\textit{offline-only stabilization} with supremizers)\\
\hline\hline
\end{tabular}\vspace{0.3cm}
\captionof{table}{Navier-Stokes problem: computational details for physical and geometrical parameters.}
\label{tab:comput_steady_NS_GEO} 
\end{center}
\section{Conclusion and perspectives}
\label{sec: Conclusion}
In this work we have used classical residual based stabilization techniques \cite{Hughes1986,Brooks1982} to develop stabilized RB methods for parametrized steady Stokes and Navier-Stokes problem. While such stabilizations have been used to handle advection dominated problems both at FE \cite{Hughes1980} and RB \cite{PR014a} level, in this work we have focused solely on the issue of inf-sup stability and pressure recovery. Indeed, construction of stable reduced basis during the offline stage does not guarantee a stable system in the online stage. The most widespread method in the RB community to tackle this issue is instead to perform a supremizer enrichment of the reduced velocity space \cite{Veroy2007}. The goal of the work has been to combine and compare the two resulting RB inf-sup stabilization options, i.e. supremizer enrichment and residual based stabilization; this gave rise to four different reduced order methods.

The first conclusion of this work is that, whenever FE stabilization is used, the \textit{offline-online stabilization} is the most appropriate way to perform the online phase. Indeed, the other option, i.e. \textit{offline-only stabilized} method, resulted in sensibly larger errors, for both velocity and pressure. We claim that this is because of the lack of consistency between the offline and online phases. Similar conclusions were reached by \cite{Torlo2018, PR014a}, although the difference is exacerbated in our case, due to the role of the pressure.

The second conclusion is on the role of supremizer enrichment. On one hand, addition of supremizers to \textit{offline-only stabilized} case does not help in improving the results, which are affected in a greater way by the lack of consistency. One the other hand, \textit{offline-online stabilization} with and without supremizers frequently showed comparable results in the error analysis, with values lower than $10^{-4}$ in all numerical examples. The case with supremizer often results in better pressure, while the case without supremizer yields better velocity in the nonlinear case. However, the differences in terms of error are often limited to one order of magnitude, still resulting in a reduced solution which would be sufficiently accurate for many practical applications. Therefore, when compared to the existing methodology based on supremizer enrichment alone, our current work allows one to make a tradeoff between accuracy and performance. Better accuracy is obtained combining residual based stabilization and supremizer enrichment, while better performance could have been obtained neglecting the enrichment step.
\par 
We still have some open questions and perspectives to improve this work in future, in order to make this approach applicable to more and more complex problems. For instance:
\begin{itemize}
\item[\textit{(i)}] in order to develop a certified stabilized RB method, an a \textit{posteriori} error analysis \cite{RB2016,Veroy2005} is needed for residual based stabilization in a reduced order setting, for which we suggest to have a look into the error analysis of stabilized FE methods \cite{Stenberg2015,Kay1999};
\item[\textit{(ii)}] the computational cost of stabilized RB method in case of nonlinear problems can be decreased by using the Empirical Interpolation method (EIM) \cite{Barrault2004};
\item[\textit{(iii)}] one can extend this work to develop a Variational MultiScale (VMS) method for turbulent flows with moderate-higher Reynolds number \cite{StabileBallari2018};
\item[\textit{(iv)}] this work is applicable to optimal control problems, see for instance environmental applications in marine sciences \cite{Maria2018};
\item[\textit{(v)}] time-dependent problems are of interest for some applications of the stabilized RB methodology;
\end{itemize}
\section*{Acknowledgements}
This work has been supported by the European Union Funding for Research and Innovation -- Horizon 2020 Program -- in the framework of European Research Council Executive Agency: H2020 ERC CoG 2015 AROMA-CFD project 681447 ``Advanced Reduced Order Methods with Applications in Computational Fluid Dynamics'' (PI Prof. Gianluigi Rozza). We also acknowledge the INDAM-GNCS project `` Advanced intrusive and non-intrusive model order reduction techniques and applications'', 2018.

\bibliographystyle{abbrv}
\bibliography{References} 
\end{document}